\documentclass[11pt]{amsart}
\usepackage{amssymb,amsthm,amsmath}
\usepackage[numbers,sort&compress]{natbib}
\usepackage{color}
\usepackage{graphicx}
\usepackage{tikz}
\usepackage[colorlinks,linkcolor=blue,anchorcolor=blue,
citecolor=blue]{hyperref}
\usepackage{comment}

\date{}

\title[cantor spectrum]{
Construction of finite differentiable quasi-periodic Schr\"odinger operators with cantor spectrum}

\author{Jiawei He}
\address{
Fujian Key Laboratory of Financial Information Processing, Putian University\\
Fujian Putian, 351100, P.R. China} \email{hermit\_well@163.com}

\author{Hongyu Cheng}
\address{
School of Mathematical Sciences, Tiangong University\\
Tianjin, 300378, P.R. China} \email{hychengmath@tiangong.edu.cn}

\newtheorem{theorem}{Theorem}[section]
\newtheorem{proposition}{Proposition}[section]
\newtheorem{lemma}{Lemma}[section]
\theoremstyle{definition}
\newtheorem{remark}{Remark}[section]
\newtheorem{claim}{Claim}[section]

\newcommand{\mi}{\mathrm{i}}     
\newcommand{\ZZ}{\mathbb{Z}}
\newcommand{\NN}{\mathbb{N}}
\newcommand{\RR}{\mathbb{R}}
\newcommand{\QQ}{\mathbb{Q}}
\newcommand{\CC}{\mathbb{C}}
\newcommand{\TT}{\mathbb{T}}
\newcommand{\mc}{\mathcal}

\numberwithin{equation}{section}
\newcommand{\red}{\textcolor{red}}
\newcommand{\blue}{\textcolor{blue}}


\begin{document}

\maketitle
\begin{abstract}
In this paper, we present a approach for
explicitly constructing quasi-periodic Schr\"odinger operators with Cantor
spectrum with $C^k$ potential. Additionally, we  provide polynomial asymptotics on the size of spectral gaps.
\end{abstract}

\section{Introduction }\label{Introduction}
Consider the quasi-periodic Schr\"odinger operator with $C^k$ potential:
\begin{equation}\label{1.1}
\left(H_{V, \alpha, \theta} x\right)_n=x_{n+1}+x_{n-1}+V(\theta+n \alpha) x_n, n \in \mathbb{Z}, d\in\NN,
\end{equation}
where $\theta \in \mathbb{T}^d$ is called the initial phase, $\alpha \in \mathbb{R}^d$ (rationally independent) is called the frequency and $V \in C^k\left(\mathbb{T}^d, \mathbb{R}\right)$ is called the potential. Let $D C_{d}(\gamma, \tau)(\gamma, \tau>0)$ denote the set of Diophantine vectors $\alpha \in \mathbb{R}^d$, i.e.,
$$
|\langle n, \alpha\rangle| \geq \frac{\gamma}{|n|^\tau}, \quad 0 \neq n \in \mathbb{Z}^d .
$$
 The Schr\"odinger operator $H_{V,\alpha,\theta}$ defined by \eqref{1.1} will be self-adjoint if we assume $V$ is real-valued.

Attributed to these  operators' close connections to physics, numerous techniques have been developed over the past forty years to study their  spectrum's structure and spectral type. We say that \eqref{1.1}  has
Cantor spectrum if their spectrum is a Cantor set. Since the 1980's, Cantor spectrum problems of Schr\"odinger operators have been the focus of significant research. The Cantor spectrum is well known to not exist when the potential is periodic. However, for almost periodic potentials, the operator could have Cantor spectrum \cite{Bellissard1982, Moser1981}, which is first conjectured by Simon \cite{SIMON1982463}.

For the quasi-periodic Schr\"odinger operators, it should be noted that if the potential in \eqref{1.1} is analytical, there has been constant progress.
In the region
of positive Lyapunov exponent, it was proved by Goldstein and Schlag \cite{Goldstein2005OnRA} that  the spectrum is a Cantor set for
any analytic potential with almost every frequency.
In the perturbative regime,
Eliasson \cite{Eliasson92} demonstrated that, for given Diophantine frequencies $\alpha \in D C(\gamma, \tau)$, $H_{V, \alpha,\theta}$ has Cantor spectrum for generic small analytic $V.$ And Puig \cite{Puig06} further extended this result to quasi-periodic Schr\"odinger operators with a fixed Diophantine frequency and nonperturbative regime potential.

The Cantor spectrum problem is more sophisticated if the potential is just finite smooth.
Cantor spectrum is a $C^0$ generic phenomenon in many classes of potentials \cite{Avilaa09,Avron81,Damanik2015}.
The work of Avila-Bochi-Damanik \cite{AvilaD08} in particular demonstrated that the Cantor spectrum holds for generic $C^0$ quasi-periodic potentials for any given rational independent frequencies.
Cai and Ge \cite{Cai2017ReducibilityOF} proved Cantor
spectrum for generic small $C^k$ quasiperiodic potential with $\alpha \in \mathrm{DC}_{d}(\gamma, \tau)$. Their arguments and the aforementioned analytical cases, however, are not constructive, are unable to offer any specific examples.

The most significant breakthroughs was made by Avila-Krikorian \cite{Avila2003ReducibilityON} and Avila-Jitomirskaya \cite{Avila2009TheTM}, who solved the Ten Martini problem \cite{Simon1983KotaniTF} associated with the
most notable almost Mathieu operators (AMO):
\begin{equation*}
\left(H_{\lambda, \alpha, \theta} x\right)_n=x_{n+1}+x_{n-1}+2\lambda\cos(\theta+n \alpha) x_n, n \in \mathbb{Z}.
\end{equation*}
They demonstrated that the AMO $H_{\lambda, \alpha, \theta}$ always possess Cantor spectrum for nonzero $\lambda$. Puig \cite{puig} previously provided the Cantor spectrum for AMO with Diophantine frequencies and small $\lambda.$

We should also mention Sinai \cite{Chulaevsky1989AndersonLF} and Wang-Zhang's work \cite{Wang2014CantorSF}, in which the Cantor spectrum was demonstrated for $H_{V, \alpha,\theta}$ in \eqref{1.1} with sufficiently large and $C^2$ cosine-like $V$. In recent work \cite{HouS19}, the authors presented a strategy for explicitly constructing a class of Schr\"odinger operators with small Gevrey quasiperiodic potential that has the Cantor spectrum.
 These are the only cases of the Cantor spectrum problem that are currently known. As far as we know, there are no  explicit example with  Cantor spectrum for $C^k$ quasi-periodic potential except for the Schr\"odinger operators with $C^2$ cosine-like $V$. Motivated by the works above, we focus our attention on the concrete $C^k$ quasi-periodic operators that have the Cantor spectrum.

\begin{theorem}\label{main2}
Let $\alpha \in \mathrm{DC}_d(\gamma, \tau)$ with $0<\gamma\leq1, \tau>d-1,$ and $V(\cdot)\in C^k\left(\mathbb{T}^d, \mathbb{R}\right)$ with $k \geq 190\tau.$ One can construct explicitly a set $\mathcal{K} \subseteq \mathbb{Z}^d$ depending on $\alpha$ and $k$ (one can refer \eqref{202207160}-\eqref{202207162} for details), such that the Schr\"odinger operator $H=H_{V, \alpha, \theta}$ with $k$ times differentiable potentials
\begin{equation}
V(\theta)=\sum_{n \in \mathcal{K}} n^{-k} \cos \langle n, \theta\rangle,
\end{equation}
has Cantor spectrum.
\end{theorem}
For the cocycle $(\alpha, A)$ with $A=e^C, C\in sl(2,\RR),$ it is uniformly hyperbolic if and only if the eigenvalues of $C$
are nonzero real numbers.
Given two cocycles $(\alpha, A_1),(\alpha, A_2) \in \mathbb{T}^d \times C^{k^{\prime}}(\mathbb{T}^d, S L(2, \mathbb{R}))$, one says that they are $C^k$ conjugated if there exists $Z \in C^k(2 \mathbb{T}^d, SL(2, \mathbb{R})),$ such that
$$
Z(\theta+\alpha) A_1(\theta) Z^{-1}(\theta)=A_2(\theta) .
$$
We say $(\alpha, A)$ is $C^{k^{\prime}, k}$ almost reducible, if $A \in C^{k^{\prime}}(\mathbb{T}^d, S L(2, \mathbb{R}))$, and the $C^k$-closure of its $C^k$ conjugacies contains a constant. We say $(\alpha, A)$ is $C^{k^{\prime}, k}$ reducible, if $A \in$ $C^{k^{\prime}}(\mathbb{T}^d, S L(2, \mathbb{R}))$ and its $C^k$ conjugacies contain a constant.
It is obvious that
a reducible system is uniformly hyperbolic if and only if the reduced constant
one is uniformly hyperbolic.

Reducibility of the cocycles is an important technique to investigate the spectral theory of quasi-periodic Schr\"odinger operators \cite{Artura10,AvilaA15,Artura23,You2019}.
When $A$ is analytic, fruitful reducibility results were obtained by the Kolmogorov-Arnold-Moser (KAM) technique, we refer the readers for \cite{Dinaburg75,Eliasson92} and the references therein. The reducibility results in lower topology are related restricted due to the notorious small divisor problem, the celebrated works are \cite{Krikorian09,Cai2017ReducibilityOF} and the references. Here, we present our quantitative version of reducibility and it will be beneficial in gap estimation as one will see later.
\begin{theorem}\label{p1}
Consider the Schr\"odinger cocycle $(\alpha, A)$ of the Schr\"odinger operator $H_{V, \alpha, \theta}$ defined in Theorem~\ref{main2}. Assume that $\alpha \in \mathrm{DC}_d(\gamma, \tau),$ $2 \rho(\alpha, A)-\langle n_J, \alpha\rangle \in \mathbb{Z}$ for $n_J \in \mathcal{K} $, and $(\alpha, A)$ is not uniformly hyperbolic, then $(\alpha, A)$ is $C^{k,k_{0}}$ reducible with  $k_0 \leq k-90\tau.$ More concretely,
there exists $B \in C^{k_{0}}(\mathbb{T}^d, \operatorname{PSL}(2, \mathbb{R}))$, such that
$$
B(\cdot+\alpha)^{-1}A(\cdot) B(\cdot)=\left(\begin{array}{cc}
1 & \zeta\\
0 & 1
\end{array}\right),
$$
with the estimates
\begin{equation}\label{est12}
\|B\|_{k_{0}}\leq n_{J}^{k_0+36\tau}, \quad n^{-(k+5\tau)}_{J}\leq|\zeta|\leq n^{-(k-62\tau)}_{J}.
\end{equation}
\end{theorem}

It is well known that the resolvent set is open, thus is an union
of open intervals, which are called spectral gaps.
For a given gap, there is a unique $n$ such that $\rho(E)=\frac{1}{2}\langle n, \omega\rangle$ holds for all $E$ in the gap \cite{Moser1981}, and $n$ is called the labeling of the gap. We denote by $\Sigma_{V, \alpha}$ and $I_n(V)=(E_n^{-}, E_n^{+})$ the spectrum and the gap with label $n$ of $H_{V, \alpha, \theta},$ respectively.
By Moser P\"oschel argument \cite{Moser1981,Moser84,zhaozhiyan},  $\zeta\neq 0 $ is a crucial prerequisite for proving that $H_{V, \alpha, \theta}$  is a Cantor set.
If $\zeta<0$, the energy of  Schr\"{o}dinger cocycle lies at left edge point of a gap	and $\zeta>0$ if  the Schr\"odinger cocycle's energy is located at the gap's right edge point. If and only if the matching energy is in a collapsed spectral gap, though, does $\zeta=0.$ Note that $\zeta\neq0 $ in \eqref{est12},  then for any $  n_J\in\mathcal{K}$, there is a open spectral gap labeled with $n_J$.   Due to  the construction of $\mathcal{K} \subseteq \mathbb{Z}^d$ (see \eqref{202207160}),
Theorem~\ref{main2} is a direct consequence of Theorem~\ref{p1}.

The precise gap estimates shed light on important physical phenomena and mathematical properties of quantum systems.
In physics, \eqref{1.1} has garnered ongoing interest since it serves as a model for the quantum Hall effect. In particular, Thouless and his coauthors \cite{Thouless8} provided a theoretic explanation of the quantization of the Hall conductance by Laughlin's argument, i.e., the Hall conductance is quantized whenever the Fermi energy lies in an energy gap, after Von Klitzing's discovery of the quantum Hall effect \cite{Klitzing80}. They did this by assuming that all gaps are open for almost Mathieu operators.  This work contributed to win Thouless the 2016 Nobel Prize.
From the perspective of mathematics, gap estimates are an essential topic in the spectrum theory of quasi-periodic Schr\"{o}dinger operators. Lower bounds estimates on spectral gaps (for all labels) are more complicated than the well-known "Dry Ten Martini Problem," whereas upper bound estimates provide a productive manner of proving spectrum homogeneity \cite{Sodin95,Sodin97}, which is a crucial subject in the study of inverse spectral theory. Our main result on spectral gaps reads as follows.

\begin{theorem}\label{gape}
Under the assumption of Theorem \ref{p1}, for any $  n_J\in\mathcal{K}$,
the spectral gap $I_{J}(V)=(E_{n_J}^{-},E_{n_J}^{+})$ has the following size:
$$n^{-\frac{11}{10}k-6\tau}_{J}\leq|I_{J}(V)|\leq n^{-\frac{9k}{10}+56\tau}_{J}.$$
\end{theorem}

What we want to emphasize is that the first result of the bounds on the size of the gaps was obtained in \cite{zhaozhiyan} with $C^{\omega}$ topology.
Latter, Cai and Wang \cite{CAI2021109035} extended the result to $C^{k}$ topology, they prove that the length of the corresponding spectral gap has a polynomial decay upper bound with respect to its label. In comparison to the previous work,  Theorem~\ref{gape} provides both the upper and lower bounds for any label $ n_J\in\mathcal{K}$.
\section{Preliminary and Construction of $\mathcal{K}$}\label{definitionpreliminary}

\subsection{Notations}
For a bounded analytic (possibly matrix valued) function $F(x)$ defined on $\mathcal{S}_{h}:=\{x=(x_{1},\cdots,x_{d})\in \CC^{d}||\Im x_{j}|<h,\forall j=1,\cdots,d\},$ let
$|F|_{h}=\sup_{x\in\mathcal{S}_{h}}\|F(x)\|$ and denote by
$C^{\omega}_{h}(\TT^{d},*)$ the set of all these *-valued functions (* will usually denote $\RR, sl(2,\RR), SL(2,\RR))$. We denote $C^{\omega}(\TT^{d},*)=\cup_{h>0}C^{\omega}_{h}(\TT^{d},*).$
Set $C^{k}(\TT^{d},*)$ to be the space of $k$ times differentiable
with continuous $k-$th derivatives functions. The norm is defined as
\begin{equation}\label{202207113}
\|F\|_{k}=\sup_{|\ell|\leq k,x\in\TT^{d}}|\partial^{\ell}F(x)|.
\end{equation}

Any integrable $\RR$-valued function $F$ on the $d-$dimensional torus has a Fourier expansion
\begin{equation}\label{202207130}
F(\theta) \sim \sum_{n \in \mathbb{Z}^d} \widehat{F}(n) e^{i\langle n, \theta\rangle}, \quad \widehat{F}(n) \triangleq \frac{1}{(2 \pi)^d} \oint_{\mathbb{T}^d} F(\theta) e^{-i\langle n, \theta\rangle} \mathrm{d} \theta .
\end{equation}
For convenience, sometimes we use the notation $\langle F\rangle=\widehat{F}(0)$.
Moreover, for any $K>0$, we define the truncating operators $\mathcal{T}_{K}$ as
\begin{equation*}
	\left(\mathcal{T}_{K} F\right)(x)=\sum_{k\in\mathbb{Z}^{d},|n|\leq K} \widehat{F}(n) e^{i\langle n,x\rangle},
\end{equation*}
and projection operator $\mathcal{R}_{K}$ as
\begin{equation*}
	\left(\mathcal{R}_{K} F\right)(x)=\sum_{n \in \mathbb{Z}^{d},|n|> K} \widehat{F}(n) e^{i\langle n,x\rangle}.
\end{equation*}

\subsection{Schr\"odinger operator and cocycle}\label{definitionoperator}
Given $A \in C^\omega\left(\mathbb{T}^d, \mathrm{SL}(2, \mathbb{R})\right)$ and $\alpha \in \mathbb{R}^d$ rationally independent, we define the quasi-periodic cocycle $(\alpha, A)$ :
$$
(\alpha, A):\left\{\begin{aligned}
\mathbb{T}^d \times \mathbb{R}^2 & \rightarrow \mathbb{T}^d \times \mathbb{R}^2 \\
(x, v) & \mapsto(x+\alpha, A(x) \cdot v)
\end{aligned}\right..
$$

Now, we turn back to Schr\"odinger operator $H_{V,\alpha,x}$ defined by \eqref{1.1}.
Note any formal solution $u=(u_{n})_{n\in\ZZ}$ of $H_{V,\alpha,x}u=Eu$ can be rewritten as
\begin{equation}\label{cocyle}
	\begin{split}
		\left(
		\begin{matrix}
			u_{n+1}
			\\
			\\
			u_{n}
		\end{matrix}
		\right)=S_{E}^{V}(x+n\alpha)\left(
		\begin{matrix}
			u_{n}
			\\
			\\
			u_{n-1}
		\end{matrix}
		\right),
	\end{split}
\end{equation}
where
\begin{equation*}
	\begin{split}
		S_{E}^{V}(x)=\left(
		\begin{matrix}
			E-V(x)\ \ \ \ -1
			\\
			\\
			1\ \ \ \ \ \ \ \ \ \ \ \ \ \ \ \  0
		\end{matrix}
		\right).
	\end{split}
\end{equation*}
We call $(\alpha,S_{E}^{V}(x))$ the Schr\"odinger cocycle.
The iterates of $(\alpha,S_{E}^{V}(\cdot))$ are of the form $(\alpha,S_{E}^{V}(\cdot))^n=(n\alpha,  S_{E,n}^{V}(\cdot)),$ where $S_{E,n}^{V}(\cdot)$
is called as the transfer matrix and defined by
\begin{equation*}
S_{E,n}^{V}(\cdot):=
\left\{\begin{array}{l l}
S_{E,V}(\cdot+(n-1)\alpha) \cdots S_{E,V}(\cdot+\alpha) S_{E,V}(\cdot),  & n\geq 0,\\[1mm]
S_{E,V}^{-1}(\cdot+n\alpha) S_{E,V}^{-1}(\cdot+(n+1)\alpha) \cdots S_{E,V}^{-1}(\cdot-\alpha), & n <0,
\end{array}\right.
\end{equation*}
then we have
\begin{equation*}
	\begin{split}
		\left(
		\begin{matrix}
			u_{n}
			\\
			\\
			u_{n-1}
		\end{matrix}
		\right)=S_{E,n}^{V}(x)\left(
		\begin{matrix}
			u_{0}
			\\
			\\
			u_{-1}
		\end{matrix}
		\right).
	\end{split}
\end{equation*}

\subsection{Rotation number} \label{rotation}
Assume that $A\in{C}^{\omega}\left(\mathbb{T}^{d},SL(2,\RR)\right)$ is homotopic to identity and introduce the map:
\begin{equation*}
       \begin{aligned}
       F:\mathbb{T}^{d}\times S^{1} \rightarrow \mathbb{T}^{d}\times S^{1},\ \
	   (x,v)\mapsto (x+\alpha,\frac{A(x)v}{\|A(x)v\|}),
       \end{aligned}
\end{equation*}
which admits a continuous lift $\widetilde{F}:\mathbb{T}^{d}\times \RR \rightarrow \mathbb{T}^{d}\times \RR$ of the form $\widetilde{F}(x,y)=(x+\alpha,y+f(x,y))$ such that $f(x,y+1)=f(x,y)$ and $\pi(y+f(x,y))=\frac{A(x)\pi(y)}{\|A(x)\pi(y)\|}$. We call that $\widetilde{F}$ is a lift for $(\alpha,A)$. Since $x\mapsto x+\alpha$ is uniquely ergodic on $\TT^{d}$, we can invoque a theorem by M.R. Herman \cite{H83} and Johnson-Moser \cite{Johnson82}: for every $(x,y)\in\TT^{d}\times\RR$ the limit
\begin{equation*}
\lim_{n\rightarrow\infty}\frac{1}{n}
\sum_{k=0}^{n-1}f(\widetilde{F}^{k}(x,y)),
\end{equation*}
exists, is independent of $(x,y)$ and the convergence is uniform in $(x,y)$; the class of this number in $\TT$ (which is independent of the chosen lift) is called the $fibered$ rotation number of $(\alpha, A)$, which is denoted by $\rho_{f}(\alpha, A)$. Moreover, the rotation number $\rho_{f}(\alpha,A)$ relates to the integrated density of states $N_{V}$ as follows:
\begin{equation}\label{pe1}
	N_{V}(E)=1-2\rho_{f}(\alpha,A).
\end{equation}

For any $C\in SL(2,\RR)$, it is immediate from the definition that
\begin{equation}\label{wanfanhou} |\rho(\alpha,A)-\rho(\alpha,C)|\leq\|A(x)-C\|_{C^{0}}^{\frac{1}{2}}.
\end{equation}
See \cite{H83} and \cite{Johnson82} for details. Besides the conclusion given by \eqref{wanfanhou}, we also have two conclusions below.
\begin{lemma}\emph{\cite{Krikorian04}}\label{rotationnumberpers}
The rotation number is invariant under the conjugation
map which is homotopic to the identity. More precisely, if
$A, B:\rightarrow SL(2, \RR)$ is continuous and homotopic to the identity, then
\begin{equation*}
\mathrm{rot}(\alpha,B(\cdot+\alpha)^{-1}A(\cdot)B(\cdot))=
\mathrm{rot}(\alpha, A).
\end{equation*}
\end{lemma}

\begin{proposition}\label{proposition2.1}
If $A:\TT^{d}\rightarrow SL(2,\RR)$ is continuous and homotopic to the identity and
$E:2\TT^{d}\rightarrow SL(2,\RR)$ is defined by
\begin{equation*}
\begin{split}
E(x)=\left(
\begin{matrix}
\cos (\pi\langle r,x\rangle )\  \ -\sin (\pi\langle r,x\rangle)\\
\\
\sin (\pi\langle r,x\rangle)\ \ \ \ \cos (\pi\langle r,x\rangle)
\end{matrix}
\right),
\end{split}
\end{equation*}
then
\begin{equation*} \rho_{f}((0,E)\circ(\alpha,A)\circ(0,E^{-1}))=\rho_{f}(\alpha,A)+\frac{\langle  r,\alpha\rangle}{2}\ \mod 1.
\end{equation*}
\end{proposition}
All concepts above for $A \in C^{\omega}\left(\mathbb{T}^d, \mathrm{SL}(2, \mathbb{R})\right)$ can be defined
 similarly
for $A \in C^{k}\left(\mathbb{T}^d, \mathrm{SL}(2, \mathbb{R})\right)$.

\subsection{Continued fraction expansion}\label{sec:cdbridge}
Let us recall some properties of the irrational number.
Given an irrational number $\alpha\in (0,1),$ we define
\begin{equation*}
	a_{0}=0,\quad \alpha_{0}=\alpha,
\end{equation*}
and inductively for $k \geq 1,$
\begin{equation*}
	a_{k}=[\alpha_{k-1}^{-1}],\quad \alpha_{k}=\alpha_{k-1}^{-1}-a_{k},
\end{equation*}
where $[\alpha]:=\max\{m\in \mathbb{Z}:\,m\leq \alpha\}$.

Let $p_{0}=0,p_{1}=1,q_{0}=1,q_{1}=a_{1},$ and recursively,
\begin{equation*}
	\begin{split}
		p_{k} = a_{k}p_{k-1}+p_{k-2},\quad
		q_{k} = a_{k}q_{k-1}+q_{k-2}.
	\end{split}
\end{equation*}
Then $ \{q_{n}\} $ is the sequence of denominators of the best rational approximations for $\alpha$. It satisfies
\begin{equation*}
	\|k\alpha\|_{\RR/\ZZ}\geq\|q_{n-1}\alpha\|_{\mathbb{T}}, \,\,\text{ for } 1 \leq k<q_{n},
\end{equation*}
and
\begin{equation}\label{202307130}
	\frac{1}{q_{n}+q_{n+1}}<\|q_{n}\alpha\|_{\RR/\ZZ}
\leq\frac{1}{q_{n+1}},
\end{equation}
where $\|x\|_{\RR/\ZZ}:=\inf\limits_{p\in\mathbb{Z}}|x-p|.$

\subsection{Construction of $\mathcal{K}$}
Set $4/5<s<1$ and let $\ell_{*}\in\NN$ with
\begin{equation}\label{202307131}
\ell_{*}=\max\{(2\|A\|)^{2.5\tau^{-1}},(5^{-1}2k)^{s^{-1}},e^{k},
(5\gamma^{-1}2^{\tau})^{\frac{\tau^{-1}}{3-2s-s^{2}}}\}.
\end{equation}
Set $M>1000$ and denote
\begin{equation}\label{202207146}
\begin{split}
\ell_{j}=M^{(1+s)^{j}},\quad j\in\NN.
\end{split}
\end{equation}
%

The set $\mathcal{K}\subset\ZZ^{d}$ is chosen such that
\begin{equation}\label{202207160}
\begin{split}
\mathrm{closure}\big\{2^{-1}\langle n,\alpha\rangle: n\in\mathcal{K}\big\}=\RR,
\end{split}
\end{equation}
\begin{equation}\label{202207161}
\begin{split}
\#\big\{n\in\mathcal{K}:\ell_{j}\leq|n|<\ell_{j+2}\big\}\leq1,
\end{split}
\end{equation}
\begin{equation}\label{2022071611}
\#\left\{n \in \mathcal{K}: 21\ell_j/10 \leq|n|<\ell_{j+1}\right\}=0,
\end{equation}
and
\begin{equation}\label{202207162}
\begin{split}
\big\{n\in\mathcal{K}:|n|<\ell_{*}\big\}=\emptyset.
\end{split}
\end{equation}

\begin{remark}
The definition of $\ell_{*}$ in \eqref{202307131} shows that $\ell_{*}$ increases with $k.$ We give this restriction to ensure that factor $\zeta$ in Theorem~\ref{p1} and the gap $|I_{J}(V)|$ in Theorem~\ref{gape} are bounded by the quantities with the parameter $k.$ Moreover, the hypothesis $M>1000$ is not optimal, neither are some other hypotheses, such as the parameters in \eqref{202307131}.
\end{remark}

In this subsection, we will construct $\mathcal{K}\subset\ZZ^{d}$ satisfying \eqref{202207160}-\eqref{202207162}. To this end, we give a auxiliary lemma to estimate the approximation of rational numbers to irrational number $\alpha\in\RR\setminus \QQ.$

\begin{lemma}\label{constructions}
Assume that $\alpha\in\RR\setminus \QQ$ with the denominators
of best rational approximations $(q_{n})_{n\in\NN}$ and
the sequence $(\ell_{j})_{j\in\NN}$ is the one defined by \eqref{202207146}, then
there exists $q_{n_{j_{*}}}\in(q_{n})_{n\in\NN}$ with
\begin{equation}\label{202306240}
\begin{split}
q_{n_{j_{*}}}\in[21\ell_{j}/20,\ 41\ell_{j}/20],
\end{split}
\end{equation}
such that
\begin{equation}\label{202306246}
\begin{split}
\|q_{n_{j_{*}}}\alpha\|_{\RR/\ZZ}:=\min_{p\in\ZZ}
|q_{n_{j_{*}}}\alpha-p|<3q_{n_j}^{-1},
\end{split}
\end{equation}
where, for the fixed $j\in\NN,$ $q_{n_{j}}$ is the one such that $q_{n_{j}}<\ell_{j}\leq q_{n_{j}+1}.$
\end{lemma}
\begin{proof}
Note that $q_{n_{j}}<\ell_{j},$ then there exists $m_{1}\in\NN$ such that
\begin{equation*}
\begin{split}
21\ell_{j}/20\leq m_{1}q_{n_{j}}<41\ell_{j}/20.
\end{split}
\end{equation*}
Moreover, set $m_{2}\in\NN$ with $m_{2}\geq m_{1},$ such that
\begin{equation*}
\begin{split}
m_{2}q_{n_{j}}\leq41\ell_{j}/20<(m_{2}+1)q_{n_{j}}.
\end{split}
\end{equation*}
The inequalities above show that
\begin{equation}\label{202306191}
\begin{split}
21\ell_{j}/20\leq m_{1}q_{n_{j}}\leq m_{2}q_{n_{j}}\leq41\ell_{j}/20,
\end{split}
\end{equation}
and
\begin{equation}\label{202306192}
\begin{split}
1<m_{1}\leq m_{2}\leq [41\ell_{j}(20q_{n_{j}})^{-1}].
\end{split}
\end{equation}
Moreover, \eqref{202307130} shows
\begin{equation*}
\begin{split}
\|q_{n_{j}}\alpha\|_{\RR/\ZZ}:=\min_{p\in\NN}|q_{n_{j}}\alpha-p|
<q_{n_{j}+1}^{-1}.
\end{split}
\end{equation*}
Set $\widetilde{q}_{j_{m}}=mq_{n}, m=m_{1},\cdots,m_{2}.$ Then
\begin{equation*}
\begin{split}
\|\widetilde{q}_{j_{m}}\alpha\|_{\RR/\ZZ}
&:=\min_{p\in\NN}|mq_{n_{j}}\alpha-p|\leq m_{2}\|q_{n_{j}}\alpha\|_{\RR/\ZZ}
\\
&<m_{2}q_{n_{j}+1}^{-1}
<3q_{n_{j}}^{-1},
\end{split}
\end{equation*}
where the last inequality is by the right inequality in \eqref{202306192} and the fact $\ell_{j}< q_{n_{j}+1}.$ Moreover, the inequalities in \eqref{202306191} show that
\begin{equation*}
\begin{split}
\widetilde{q}_{j_{m}}= mq_{n_{j}}\in[21\ell_{j}/20,\ 41\ell_{j}/20], m=m_{1},\cdots,m=m_{2}.
\end{split}
\end{equation*}
The discussions above show that $\widetilde{q}_{j_{m}}, m=m_{1},\cdots,m_{2},$ are the wanted ones.
\end{proof}

\begin{lemma}\label{construction}
For any $\alpha\in\RR^{d},$ there exists $\mathcal{K}\subset\ZZ^{d}$ satisfying \eqref{202207160}-\eqref{202207162}.
\end{lemma}
\begin{proof}
We order all $n=(n_{1},\cdots,n_{d})\in\ZZ^{d}$ via the lexicographic order $(|n|,n_{1},n_{2},\cdots,n_{d}),$
for example
\begin{equation*}
\begin{split}
n^{(0)}&=(0,\cdots,0),\ n^{(1)}=(-1,0,0,\cdots,0),\ n^{(2)}=(0,-1,0,\cdots,0),\ \cdots,\\
n^{(d+1)}&=(1,0,0,\cdots,0),\ n^{(d+2)}=(0,1,0,\cdots,0),\ \cdots.
\end{split}
\end{equation*}
Obviously,
\begin{equation}\label{202207164}
\begin{split}
(n^{(m)})_{m\in\NN}=\ZZ^{d}, \ \
|n^{(m)}|\leq m,\ m=0,1,\cdots.
\end{split}
\end{equation}

For $\alpha=(\alpha_{1},\cdots,\alpha_{d})\in DC(\gamma,\tau),$
we know that
$\alpha_{j}\in (0,1)\setminus\QQ,\forall j=1,\cdots,d.$ Without loss of generality, we fix $\alpha_{1},$ and
let $(q_{n})_{n\in\NN}$ be the sequence of denominators of the best rational approximations for $\alpha_{1}.$ Set
\begin{equation}\label{202207170}
\begin{split}
\widehat{n}^{(j)}=(q_{n_{j_{*}}},0\cdots,0),\quad j\geq1,
\end{split}
\end{equation}
where, for the fixed $j\in\NN,$ $q_{n_{j_{*}}}$ is the one constructed in Lemma~\ref{constructions} with $\alpha_{1}$ in place of $\alpha.$

Set
\begin{equation}\label{202207168}
\begin{split}
\widetilde{n}^{(m)}=n^{(m)}+\widehat{n}^{(j_{m})}, m\in\NN
\end{split}
\end{equation}
with
\begin{equation}\label{202207178}
\begin{split}
j_{1}\geq \ell_{*},\ j_{m+1}-j_{m}\geq2.
\end{split}
\end{equation}
Then, \eqref{202207164}, \eqref{202207170} and \eqref{202306240} imply that
\begin{equation*}
\begin{split}
\ell_{j_{m}}<21\ell_{j_{m}}/20-m\leq|\widetilde{n}^{(m)}|
\leq 41\ell_{j_{m}}/20+m<21\ell_{j_{m}}/10,\ m\in\NN.
\end{split}
\end{equation*}
That is
\begin{equation}\label{202206160}
\begin{split}
|\widetilde{n}^{(m)}|\in[\ell_{j_{m}},21\ell_{j_{m}}/10),\ m\in\NN.
\end{split}
\end{equation}

Set
\begin{equation*}
\begin{split}
\mathcal{K}=\{\widetilde{n}^{(m)}:  m\in\NN\}.
\end{split}
\end{equation*}

Now we verify that the set $\mathcal{K}$ satisfies the estimates \eqref{202207160}-\eqref{202207162}. The relations in \eqref{202207178} and \eqref{202206160} yield \eqref{202207161}-\eqref{202207162}.
Moreover, \eqref{202207168} yields
\begin{equation}\label{202207179}
\begin{split}
\langle\widetilde{n}^{(m)},\alpha\rangle= \langle n^{(m)},\alpha\rangle+q_{n_{j_{m_{*}}}}\alpha_{1},\ \forall m\in\NN.
\end{split}
\end{equation}
Set $q_{n_{j_{m}}}<\ell_{j_{m}}\leq q_{n_{j_{m}}+1},$ then \eqref{202306240} shows that
\begin{equation*}
\begin{split}
\lim_{m\rightarrow\infty}\|q_{n_{j_{m*}}}\alpha_{1}\|_{\RR\setminus\QQ}&
\leq\lim_{n\rightarrow\infty}2q_{n_{j_{m}}}^{-1}=0,
\end{split}
\end{equation*}
which, together with \eqref{202207179} and \eqref{202207164}: $(n^{(m)})_{m\in\NN}=\ZZ^{d},$ imply
\begin{equation*}
\begin{split}
\mathrm{closure}\big\{2^{-1}\langle \widetilde{n}^{(m)},\alpha\rangle: \widetilde{n}^{(m)}\in\mathcal{K}\big\}=\mathrm{closure}\big\{2^{-1}\langle n^{(m)},\alpha\rangle\}_{m\in\ZZ^{d}}=\RR,
\end{split}
\end{equation*}
which yields \eqref{202207160}.

\end{proof}

\section{Finite almost reducibility}\label{almostreducibility}
In this section, we will establish our main KAM induction and then give the basic quantitative estimates in the case of reducibility and almost reducibility, respectively. These estimates will be applied to control the growth of corresponding Schr\"odinger cocycles.

\subsection{Auxiliary lemmas}
Denote by $M_{n}(\CC)$ ($M_{n}(\RR)$) the set of complex (real) $n$ by $n$ matrices and and by $GL(n,\CC)$ ($GL(n,\RR)$) the set of matrices in $M_{n}(\CC)$ ($M_{n}(\RR)$) with nonvanishing determinant.
The groups $SL(2,\RR)$ and $SU(1,1),$ the subgroups of $M_{2}(\CC),$ are isomorphic by $M=\frac{1}{1+\mi}\left(
\begin{matrix}
1\  -\mi
\\
1\ \  \ \ \mi
\end{matrix}
\right),$ that is $SU(1,1)=MSL(2,\RR)M^{-1}.$  Correspondingly, $su(1,1)=M sl(2,\RR)M^{-1}.$ Moreover, set $PSL(2,\RR)$ and $PSU(1,1),$ by the quotient groups defined by
\begin{equation*}
\begin{split}
PSL(2,\RR):=SL(2,\RR)/\{\pm I\},\quad PSU(1,1):=SU(1,1)/\{\pm I\}.
\end{split}
\end{equation*}
%

\begin{lemma}\label{sS}
 For   $C\in su(1,1)$ with $spec (C)=\{\pm\lambda\}$, and  \begin{equation}
C=\left(\begin{array}{cc}
i a & b \\
\overline{b} & -i a
\end{array}\right), \quad \text{where }  a \in \mathbb{R}, \quad b \in \mathbb{C},
\end{equation}
 then $A\triangleq e^C \in SU(1,1) $ and with the following form
 \begin{equation}\label{202305300}
A=\left(\begin{array}{cc}
\cosh (\lambda)+i a \frac{\sinh (\lambda)}{\lambda} & b \frac{\sinh (\lambda)}{\lambda} \\
\overline{b} \frac{\sinh (\lambda)}{\lambda} & \cosh (\lambda)-i a \frac{\sinh (\lambda)}{\lambda}
\end{array}\right).
\end{equation}
\end{lemma}
\begin{proof}
    By simple calculations, it follows.
\end{proof}
\begin{lemma}\label{RD}
Let $A=\left(\begin{array}{cc} a& {b} \\ \overline{b} & \overline{a}\end{array}\right)\in SU(1,1)$ with $\operatorname{spec}(A)=\{ 1\},$ where $|a|^2-|b|^2=1$ and $  a, b \in \mathbb{C},$  then there exists $\phi\in \TT$, $R_{\phi}:=\left(\begin{array}{cc}
\cos 2 \pi \phi & -\sin 2 \pi \phi \\
\sin 2 \pi \phi & \cos 2 \pi \phi
\end{array}\right)$, satisfying $R_{-\phi}M^{-1} A MR_{\phi}=\left(\begin{array}{cc}1 & |b| \\ 0 & 1\end{array}\right)$.
\end{lemma}

\begin{proof}
Note that $\operatorname{spec}(A)=\{ 1\},$  $Re(a)=1$. Combing with $|a|^2-|b|^2=1$, we have $Im(a)=|b|$. Thus
$A=\left(\begin{array}{cc} 1+i|b|& {b} \\ \overline{b} & 1-i|b|\end{array}\right)$,
the desired result can be reached with some straightforward computations.
\end{proof}

\subsection{Estimates on Algebraic Conjugations}

For a given $B \in \mathfrak{B}_h\left(\mathbb{T}^d, P S U(1,1)\right)$, we define the operator $A d(B)$ as
$$
A d(B) . W \triangleq B W B^{-1}, \quad W \in \mathfrak{B}_h(\mathbb{T}^d, s u(1,1)) .
$$
Let
\begin{equation}\label{312}
  W=\left(\begin{array}{cc}
i u & w \\
\bar{w} & -i u
\end{array}\right) \in \mathfrak{B}_h(\mathbb{T}^d, s u(1,1)),
\end{equation}
and write $B W B^{-1}$ as
\begin{equation}\label{313}
  W_{+} \triangleq B W B^{-1}=\left(\begin{array}{cc}
i u_{+} & w_{+} \\
\bar{w}_{+} & -i u_{+}
\end{array}\right) .
\end{equation}
In the following, we give estimates of $u_{+}, w_{+}$.

\begin{lemma}\label{3.4}
 Let $A \in SU(1,1)$ with $\operatorname{spec}(A)=\{ e^{\pm i \rho}\}, \text{where } \rho\neq 0
 , P \in S U(1,1)$ satisfying $P A P^{-1}=\left(\begin{array}{cc}e^{i \rho} & 0 \\ 0 & e^{-i \rho}\end{array}\right)$. Write $W \in \mathfrak{B}_h\left(\mathbb{T}^d, s u(1,1)\right)$ in the form
\eqref{312} and $W_{+}=P W P^{-1}$ in the form \eqref{313}. Then
\begin{equation}\label{3.31}
   \left|\left\langle w_{+}\right\rangle\right|  \geq \frac{1}{2}|\rho|^{-1}\|P\|^{-2}\|[A,\langle W\rangle]\|,
\end{equation}
\begin{equation}\label{3.32}
\begin{aligned}
\mid \widehat{w}_{+}(n) \mid &\geq \frac{\|P\|^2+1}{2}(|\widehat{w}(n)|-3 \max \{|\widehat{w}(-n)|,|\widehat{u}(n)|\}) \\
& \geq|\widehat{w}(n)|-3 \max \{|\widehat{w}(-n)|,|\widehat{u}(n)|\}, \ \ \forall n \in \mathbb{Z}^d.
\end{aligned}
\end{equation}

\end{lemma}
\begin{proof}
The proof is essentially contained in Lemma 3.4 of \cite{HouS19}; however, it should be noted that in Lemma 3.4 of \cite{HouS19}, $A \in su(1,1)$, whereas here it also holds for $A \in SU(1,1)$.
\end{proof}

\begin{lemma}\label{3.5}\cite{HouS19}
 Let $W \in s u(1,1)$ and $B \in \mathfrak{B}_h\left(\mathbb{T}^d, P S U(1,1)\right)$ satisfying $|B-I|_h \leq \frac{1}{2}$. Write $W_{+}=B^{-1} W B$ in the form \eqref{313}. Then,
$$
|\widehat{w}_{+}(n)| \geq|\widehat{w}(n)|-4|B-I|_h|W|_h e^{-|n| h}, \ \ \forall n \in \mathbb{Z}^d.
$$
\end{lemma}

\subsection{Normal form}
The quasi-periodic cocycles defined in \eqref{cocyle} can be rewritten as
\begin{equation}
\left(\begin{array}{c}
u_{n+1} \\
u_n
\end{array}\right)=\left(A_0+F_0(\theta+n \alpha)\right)\left(\begin{array}{c}
u_n \\
u_{n-1}
\end{array}\right)
\end{equation}
with
$$
A_0=
\left(\begin{array}{cc}
E & -1 \\
1 & 0
\end{array}\right),\ F_0=
\left(\begin{array}{cc}
-V & 0 \\
0 & 0
\end{array}\right).
$$
Furthermore, take $\left\{n_j\right\}_{j \in \mathbb{N}} \in \mathcal{K} \subset \mathbb{Z}^d$
$$
V(x)=\sum_{j \in \mathbb{N}} |n_j|^{-k} \cos (2 \pi\langle n_j, x\rangle).
$$
If we denote
\begin{equation}\label{111}
\widetilde{W}=
\left(\begin{array}{cc}
0 & 0 \\
1 & 0
\end{array}\right),\ W=M\widetilde{W}M^{-1},
\end{equation}
then we also have
\begin{equation}\label{11}
S_E^V(x):=\left(\begin{array}{cc}
E-V(x) & -1 \\
1 & 0
\end{array}\right)
=A_0+F_0=A_{0}e^{\widetilde{F}}=A_{0} \prod_{n_j \in \mathcal{K}}e^{V_{j}(x)\widetilde{W}},
\end{equation}
where
\begin{equation*}
V_{j}(x)= |n_j|^{-k} \cos (2 \pi\langle n_j, x\rangle),\quad \widetilde{F}=
\left(\begin{array}{cc}
0 & 0 \\
V & 0
\end{array}\right).
\end{equation*}
\begin{remark}
Note that for the potential $V(x),$ we set its Fourier coefficients as $|n_j|^{-k}, n_j\in \mathcal{K}.$ Actually, the Fourier coefficients are not fixed and can be any real numbers bounded by $c|n_j|^{-k}, n_j\in \mathcal{K}.$ Thus, we have constructed a family of potentials such that the corresponding Schr\"odinger operators possess Cantor spectrum.
\end{remark}

\subsection{One Step of KAM}

In this subsection, our main aim  is to find a conjugation to simplify the cocycle $(\alpha, Ae^F)$.
Assume that for given $\eta>0,\ \alpha\in\RR^{d}$ and $A\in SU(1,1),$  we define a decomposition
$C_{h}^{\omega}(\TT,su(1,1))=C_{h}^{nre}(\TT,su(1,1))\oplus C_{h}^{re}(\eta)$ satisfying that for any $Y\in C_{h}^{nre}(\eta),$
\begin{equation}\label{202207240}
\begin{split}		
A^{-1}Y_{j}(\cdot+\alpha)A-Y_{j}(\cdot)\in C_{h}^{nre}(\eta),
|A^{-1}Y_{j}(\cdot+\alpha)A-Y_{j}(\cdot)|_{h}\geq \eta|Y|_{h}.
\end{split}
\end{equation}	
Set $\mathbb{P}_{nre}, \mathbb{P}_{re}$ be the standard projections from $C_{h}^{nre}(\TT,su(1,1))$ onto $C_{h}^{nre}(\eta)$ and
$C_{h}^{nre}(\eta),$ respectively.
Then we have the following crucial lemma which helps us remove all the non-resonant
terms:

\begin{lemma}\label{implicittheorem}
Assume that $A\in SU(1,1)$
and $\eta\leq(2\|A\|)^{-30}.$ Then, for any $F\in C^{\omega}_{h}(\TT,su(1,1))$ with $|F|_{h}<\eta^{61/30},$ there exists $Y\in C^{nre}_{h}(\eta)$, $F_{*}\in C^{re}_{h}(\eta)$ such that
\begin{equation}\label{202207226}
e^{Y(x+\alpha)}(Ae^{F(x)})e^{-Y(x)}=Ae^{F_{*}(x)},
\end{equation}
with estimates
\begin{equation}\label{202207227}
|Y|_{h}\leq2\eta^{-1}|F|_{h},\quad |F_{*}-\mathbb{P}_{re}F|_{h}\leq2\eta^{-6.1}|F|_{h}^{2}.
\end{equation}
\end{lemma}
\begin{proof}
We give the proof by Newton iteration. Decompose $F$ as
\begin{equation*}
F=:F_{0}=\mathbb{P}_{re}F_{0}+\mathbb{P}_{nre}F_{0}
=F_{0}^{(re)}+F_{0}^{(nre)}.
\end{equation*}
Assume that we have constructed $\{Y_{\ell}\}_{\ell=0}^{j-1}\subset C_{h}^{nre}(\TT,su(1,1))$ and
$\{F_{\ell}\}_{\ell=1}^{j}\subset C_{h}^{\omega}(\TT,su(1,1))$
with the estimates
\begin{equation}\label{202207228}
\begin{split}
 |Y_{p-1}|_{h}&\leq\eta^{-1}|F_{p-1}^{(nre)}|_{h},\quad |F_{p}^{(nre)}|_{h}\leq\eta^{-61/30}|F_{p-1}^{(nre)}|_{h}^{2},\\
 |F_{p}^{(re)}&-F_{p-1}^{(re)}|_{h}
 \leq\eta^{-61/30}|F_{p-1}^{(nre)}|_{h}^{2},\ p=0,\cdots,j-1,
\end{split}
\end{equation}
such that
\begin{equation}\label{202207237}
\begin{split}
e^{Y_{p}(x+\alpha)}(Ae^{F_{p}(x)})e^{-Y_{p}(x)}
=Ae^{F_{p+1}(x)},\ p=0,\cdots,j-1.
\end{split}
\end{equation}

In the following, we will construct $Y_{j}\in C_{h}^{nre}(\TT,su(1,1)),$ and $F_{j+1}\in C_{h}^{\omega}(\TT,su(1,1)),$ such that \eqref{202207228} and \eqref{202207237} hold with $p=j.$ Assume that $Y_{j}\in C^{nre}_{h}(\TT,su(1,1))$ is the solution to
\begin{equation}\label{202207230}
\begin{split}
A^{-1}Y_{j}(x+\alpha)A-Y_{j}(x)=-F_{j}^{(nre)}.
\end{split}
\end{equation}
For $Y_{j}\in C^{nre}_{h}(\eta),$ by \eqref{202207240} we know that
\begin{equation*}
\begin{split}
|A^{-1}Y_{j}(x+\alpha)A-Y_{j}(x)|_{h}
\geq\eta|Y_{j}|_{h}.
\end{split}
\end{equation*}
Thus, there is unique solution to \eqref{202207230} with estimate
\begin{equation}\label{202207231}
\begin{split}
|Y_{j}|_{h}\leq \eta^{-1}|F_{j}^{(nre)}|_{h}.
\end{split}
\end{equation}
Moreover,
\begin{equation}\label{202207232}
\begin{split}
e^{Y_{j}(x+\alpha)}&(Ae^{F_{j}(x)})e^{-Y_{j}(x)}\\
&=A\{I+F_{j}^{(re)}(x)+A^{-1}O(|Y_{j}|^{2},|F_{j}|^{2},|Y||F_{j}|)\}\\
&=Ae^{F_{j+1}(x)},
\end{split}
\end{equation}
where
\begin{equation}\label{202207233}
\begin{split}
F_{j+1}(x)&=F_{j}^{(re)}(x)
+A^{-1}O(|Y_{j}|^{2},|F_{j}|^{2},|Y||F_{j}|)\}\\
&+F_{j}^{(re)}(x)^{2}+
A^{-1}F_{j}^{(re)}O(|Y_{j}|^{2},|F_{j}|^{2},|Y_{j}||F_{j}|).
\end{split}
\end{equation}

Now, we will verify the last two inequalities in \eqref{202207228} hold with $p=j.$ Obviously, \eqref{202207231}, \eqref{202207233}, together with the fact $\eta\leq(2\|A\|)^{-30},$ imply
\begin{equation*}
\begin{split}
|F_{j+1}-F_{j}^{(re)}|_{h}\leq 2\|A^{-1}\||Y_{j}|_{h}^{2}
\leq\eta^{-1/30}\eta^{-2}|F_{j}^{(nre)}|_{h}^{2}=
\eta^{-61/30}|F_{j}^{(nre)}|_{h}^{2},
\end{split}
\end{equation*}
which yields
\begin{equation*}
\begin{split}
|F_{j+1}^{(re)}-F_{j}^{(re)}|_{h}&\leq|F_{j+1}-F_{j}^{(re)}|_{h}
\leq\eta^{-61/30}|F_{j}^{(nre)}|_{h}^{2},\\
|F_{j+1}^{(nre)}|_{h}=|(F_{j+1}-F_{j}^{(re)}&)^{nre}|_{h}
\leq |F_{j+1}-F_{j}^{(re)}|_{h}\leq\eta^{-61/30}|F_{j}^{(nre)}|_{h}^{2}.
\end{split}
\end{equation*}

Repeat the above process for infinite times, we know that \eqref{202207228} and \eqref{202207237} hold for all $p\in\NN.$
Set
\begin{equation}\label{202207234}
\begin{split}
e^{Y}:=\lim_{p\rightarrow\infty}e^{Y_{p-1}}
\cdots e^{Y_{1}}e^{Y_{0}},\ \
F_{*}=\lim_{p\rightarrow\infty}F_{p}.
\end{split}
\end{equation}
Thus, \eqref{202207237} implies
\begin{equation*}
\begin{split}
e^{Y(x+\alpha)}(Ae^{F(x)})e^{-Y(x)}
=Ae^{F_{*}(x)}.
\end{split}
\end{equation*}

In the following, we will verify that the functions $Y$ and $F_{*}$
defined above satisfy estimates in \eqref{202207227} with  $F_{*}=F_{*}^{(re)}\in C_{h}^{re}(\TT,su(1,1)).$

The second estimate in \eqref{202207228} yields
\begin{equation}\label{202207235}
\begin{split}
|F_{p}^{(nre)}|_{h}&\leq\eta^{-61/30}|F_{p-1}^{(nre)}|_{h}^{2}
\leq\eta^{-(1+2)61/30}|F_{p-2}^{(nre)}|_{h}^{2^{2}}\leq\cdots\\
&\leq
\eta^{-(1+2+\cdots+2^{p-1})61/30}|F_{0}^{(nre)}|_{h}^{2^{p}}
\leq\{\eta^{-61/30}|F_{0}^{(nre)}|_{h}\}^{2^{p}}.
\end{split}
\end{equation}
Thus,
\begin{equation*}
\begin{split}
|F_{*}^{(nre)}|_{h}=\lim_{p\rightarrow\infty}|F_{p}^{(nre)}|_{h}
\leq\lim_{p\rightarrow\infty}
\{\eta^{-61/30}|F_{0}^{(nre)}|_{h}\}^{2^{p}}=0,
\end{split}
\end{equation*}
since $|F_{0}|_{h}<\eta^{61/30}.$ That is $F_{*}=F_{*}^{(re)}\in C_{h}^{re}(\TT,su(1,1)).$

On the other hand, note
\begin{equation*}
\begin{split}
F_{p}^{(re)}-F_{0}^{(re)}=\sum_{\ell=0}^{p-1}\{F_{\ell+1}^{(re)}
-F_{\ell}^{(re)}\},
\end{split}
\end{equation*}
then, by the third inequality in \eqref{202207228} and \eqref{202207235} we get
\begin{equation*}
\begin{split}
|F_{p}^{(re)}-F_{0}^{(re)}|_{h}&
\leq\sum_{\ell=0}^{p-1}|F_{\ell+1}^{(re)}
-F_{\ell}^{(re)}|_{h}
\leq\sum_{\ell=0}^{p-1}\eta^{-61/30}|F_{\ell}^{(nre)}|_{h}^{2}\\
&\leq\sum_{\ell=0}^{p-1}\eta^{-61/30}
\{\eta^{-61/30}|F_{0}^{(nre)}|_{h}\}^{2^{\ell+1}}.
\end{split}
\end{equation*}
Thus,
\begin{equation*}
\begin{split}
|F_{*}^{(re)}-F_{0}^{(re)}|_{h}&=
\lim_{p\rightarrow\infty}|F_{p}^{(re)}-F_{0}^{(re)}|_{h}\\
&\leq\lim_{p\rightarrow\infty}\sum_{\ell=0}^{p-1}\eta^{-61/30}
\{\eta^{-61/30}|F_{0}^{(nre)}|_{h}\}^{2^{\ell+1}}\\
&\leq2\eta^{-6.1}|F_{0}^{(nre)}|_{h}^{2}.
\end{split}
\end{equation*}
The inequality above, together with the fact $F_{*}=F_{*}^{(re)},$ implies that $F_{*}$ defined by \eqref{202207234} is the wanted one. With the similar calculations above, we also get the estimates about $Y,$ we omit the details.
\end{proof}

Let us introduce some notations.
For $m\in \ZZ^{+},$ we define
\begin{equation}\label{202207121}
\begin{split}
\varepsilon_{m}=m^{-1},\ \ N_{\ell_{j}}=2\ell_{j+1},\ \ h_{j}=10\tau\ell_{j}^{-1}\ln{\ell_j}.
\end{split}
\end{equation}
Denote
\begin{equation}
\begin{aligned}
& \mathcal{Z}_j \triangleq\{n \in \mathbb{Z}^d|\ \ \ell_j \leq| n \mid<\ell_{j+1}\},  j\geq1, \\
& \mathcal{Z}_0 \triangleq\{n \in \mathbb{Z}^d|\ \ 0 \leq| n \mid<\ell_1\}.
\end{aligned}
\end{equation}

\begin{equation}\label{2022071471}
\tilde{h}_j=\left\{\begin{array}{l}
h_j, \quad \text { as } \mathcal{K} \cap \mathcal{Z}_j \neq \emptyset \\
\frac{3}{4} h_{j-1}, \quad \text { as } \mathcal{K} \cap \mathcal{Z}_j=\emptyset
\end{array}\right. .
\end{equation}

Within the above concepts,
our main result is the following:
\begin{lemma}\label{nonresonantcase}
Let $\alpha\in DC_{d}(\gamma,\tau)$, $\gamma>0,\tau>1$. Consider the cocycle $(\alpha, A_{j}e^{F_{j}(x)})$ where $A_{j}\in
 SU(1,1)$ and
 \begin{equation}\label{202207148}
 e^{F_{j}(x)}=e^{f_{j}}\prod_{p\geq j}e^{Ad(B^{(j)})(V_{p}W)}
 \end{equation}
where $W$ is the one defined in \eqref{111} and $f_{j} \in \mathfrak{B}_{h_j}(\mathbb{T}^d, su(1,1))$ satisfying
\begin{equation}\label{202201180}
\begin{split}
|f_{j}|_{\frac{3}{4}h_{j-1}}\leq \varepsilon^{2k}_{\ell_{j}},\ \
|B^{(j)}|_{h_{j}}
\leq\varepsilon_{\ell_{j}}^{-18\tau}.
\end{split}
\end{equation}
Then there exists $ B_j \in C_{h_{j}}^{\omega}(2\TT^{d},SU(1,1))$, $F_{j+1}\in su(1,1),$ and $A_{j+1}\in SU(1,1)$, such that
\begin{equation}\label{202201181}
{B}_{j}(x+\alpha)(A_{j}e^{F_{j}(x)}){B}_{j}(x)^{-1}
=A_{j+1}e^{F_{j+1}(x)},
\end{equation}
where $F_{j+1}(x)$ is the one defined by \eqref{202207148} with the estimate \eqref{202201180} with $j+1$ in place of $j.$ Moreover, the following conclusions also hold:

\textbf{Case one:} $A_{j}\in \mathcal{NR}(N_{\ell_{j}},\varepsilon^{4\tau}_{\ell_{j}}):$ we have the following estimates
\begin{equation}\label{2j}
|B_j-I|_{h_j} \leq 6\varepsilon^{k-69\tau}_{\ell_j}, \quad\|A_{j+1}-A_j\| \leq 20\varepsilon^{k-57\tau}_{\ell_{j}}\|A_j\|,
\end{equation}
\begin{equation}\label{202306020}
\begin{split}
|B^{(j+1)}|_{h_{j}}\leq2|B^{(j)}|_{h_j}\leq2\varepsilon^{-18\tau}_{\ell_{j}},
\end{split}
\end{equation}
and,
\begin{equation}\label{3j1}
\|A_j+\big\langle A d\big(B^{(j)}\big) \cdot\big(V_j W\big)\big\rangle-A_{j+1}\|\leq 80 \varepsilon_{\ell_j}^{2k-188\tau}, \ \text{as}\ \ \mathcal{K} \cap \mathcal{Z}_j \neq \emptyset,
\end{equation}

\begin{equation}\label{3j2}
\|A_{j+1}-A_j\|\leq80 \varepsilon_{\ell_j}^{2k-188\tau}, \quad \text{as}\ \ \mathcal{K} \cap \mathcal{Z}_j=\emptyset.
\end{equation}

\textbf{Case two:} $A_{j}\in \mathcal{RS}(N_{\ell_{j}},\varepsilon^{4\tau}_{\ell_{j}}),$ that is there exists  $n^{*}_{j}$ with $0<|n^{*}_{j}| \leq N_{\ell_{j}},$ such that
$|2 \rho_{j}-\langle n^{*}_{j}, \alpha\rangle|<\varepsilon^{4\tau}_{\ell_{j}}.$
Then, $B_j:=Q_{n^{*}_{j}} \breve{B}_{j} P_j,$ where
$P_j \in S U(1,1)$ such that $P_j^{-1} A_j P_j=\left(\begin{array}{cc}e^{i \rho_j} & 0 \\ 0 & e^{-i \rho_j}\end{array}\right)$ with $0<\rho_j <2\pi,$
$Q_{n^{*}_{j}}\in\mathfrak{B}_{h_j}(2\mathbb{T}^d, P S U(1,1))$
and $\breve{B}_{j}\in\mathfrak{B}_{h_j}(\mathbb{T}^d, SU(1,1))$ with the following estimates
\begin{equation}\label{6j}
| \breve{B}_{j_0} -I|_{h_j} \leq18\varepsilon_{\ell_j}^{k-65\tau}, \ \ \|P_j\|^2 \leq 3 \varepsilon^{-4\tau}_{\ell_{j}}, \ \ |Q_{n^{*}_j}|_{h_{j+1}}\leq\varepsilon^{-10\tau}_{\ell_{j+1}}.
\end{equation}
Moreover,
\begin{equation}\label{6j2}
A_{j+1} \in \mathcal{N} \mathcal{R}(N_{\ell_{j+1}}, \varepsilon^{4\tau}_{\ell_{j+1}}).
\end{equation}
If we write $A_{j+1}=\left(\begin{array}{cc}*& {b_{j+1}} \\ \overline{b_{j+1}} & *\end{array}\right)$ and $P_j A d\left(B^{(j)}\right) \cdot\left(v_j W\right) P_j^{-1}=\left(\begin{array}{cc}* & T_j \\ \overline{T_j} & *\end{array}\right)$, we furthermore have
\begin{equation}\label{7j}
\left\{\begin{array}{l}
\big|b_{j+1}-\widehat{T}_j\big(n^{*}_j\big)\big| \leq 243 \varepsilon_{\ell_j}^{2k-147\tau} e^{-|n^{*}_j| h_j}, \quad \mathcal{K} \cap \mathcal{Z}_j \neq \emptyset; \\
\big|b_{j+1}\big| \leq 243 \varepsilon_{\ell_j}^{2k-147\tau} e^{-\frac{3}{4}|n^{*}_j| h_{j-1}}, \quad \mathcal{K} \cap \mathcal{Z}_j=\emptyset .
\end{array}\right.
\end{equation}

\end{lemma}

\begin{proof}
Define
\begin{equation}\label{202207150}
\begin{split}
\Lambda=\{f\in C_{h_{j}}^{\omega}(\TT^{d},su(1,1))|
~f(x)=\sum_{n\in \mathbb{Z}^{d},0<|n|\leq N_{\ell_{j}}}\widehat{f}(n)e^{i\langle k,x\rangle}\},
\end{split}
\end{equation}
and set
\begin{equation}\label{202207155}
\begin{split}
e^{\Tilde{F}_{j}(x)}=e^{f_{j}(x)}e^{Ad(B^{(j)})V_{j}(x)W},
\end{split}
\end{equation}
where
\begin{equation}
V_j=\left\{\begin{array}{l}
n^{-k}_{j}\cos \left(2 \pi\left\langle n_{j}, x\right\rangle\right), \ \text { if } \{n_{j}\} \triangleq \mathcal{K} \cap \mathcal{Z}_j,  \\
0, \quad \text { if } \mathcal{K} \cap \mathcal{Z}_j=\emptyset.
\end{array}\right.
\end{equation}
Now, we prove that $\Tilde{F}_{j}\in C_{\tilde{h}_{j}}^{\omega}(\TT^{d},su(1,1))$ with certain estimates. First, we consider $\mathcal{K} \cap \mathcal{Z}_j \neq \emptyset$, then \eqref{2022071611} implies
\begin{equation}\label{v1}
\begin{split}
|V_{j}(x)|_{\tilde{h}_j}
&\leq|\widehat{V}(n_{j})|e^{| n_{j}h_{j}|}
=|n_{j}|^{-k}e^{| n_{j}h_{j}|}\\
&
=\exp\{-k\ln|n_{j}|+10\tau \ell_{j}^{-1}|n_{j}|\ln\ell_{j}\}\\
&\leq\exp\{-k\ln\ell_{j}+21\tau\ln\ell_{j}\}\\
&= \ell_{j}^{-(k-21\tau)}\leq\varepsilon_{\ell_{j}}^{k-21\tau}.
\end{split}
\end{equation}
Then, the inequality above, together with the estimates in \eqref{202201180}, yields
\begin{equation*}\label{202207156}
\begin{split}
|\Tilde{F}_{j}|_{\tilde{h}_j}&\leq
2(|f_{j}|_{\frac{3}{4}h_{j-1}}+|B^{(j)}|_{h_{j}}^{2}
|V_{j}|_{h_{j}})\\
&\leq 2(\varepsilon^{2k}_{\ell_{j}}+\varepsilon_{\ell_{j}}^{(k-21\tau)-36\tau})
\leq3\varepsilon^{k-57\tau}_{\ell_{j}}.
\end{split}
\end{equation*}
For the case $\mathcal{K} \cap \mathcal{Z}_j = \emptyset$, $\Tilde{F}_{j}=f_j$.   By \eqref{202201180}, we also have $
|\Tilde{F}_{j}|_{\tilde{h}_j}
\leq3\varepsilon^{k-57\tau}_{\ell_{j}}.
$
Thus in both case, we have
\begin{equation}\label{estimate_F}
|\Tilde{F}_{j}|_{\tilde{h}_j}
\leq3\varepsilon^{k-57\tau}_{\ell_{j}}.
\end{equation}

\textbf{Case one:} $A_{j}\in \mathcal{NR}(N_{\ell_{j}},\varepsilon^{4\tau}_{\ell_{j}})$,
 for any $Y_j\in\Lambda_{j},$
\begin{equation}\label{202207241}
\begin{split}		
A_j^{-1}Y_{j}(x+\alpha)A_j-Y_{j}(x)\in C_{{\tilde{h}_j}}^{nre}(\varepsilon^{4\tau}_{\ell_{j}}).
\end{split}
\end{equation}	
Note $\alpha\in DC_{d}(\gamma,\tau),$  we have
\begin{equation}\label{202207151}
\begin{split}
|\langle n,\alpha\rangle|\geq\gamma|n|^{-\tau}\geq\varepsilon_{\ell_{j}}^{4\tau},  0<|n|<N_{\ell_{j}}.
\end{split}
\end{equation}
Moreover, since we assume
$A_j\in\mathcal{NR}(N_{\ell_{j}},
\varepsilon_{\ell_{j}}^{4\tau}),$ then
\begin{equation*}
\begin{split}
|2\rho_{j}\pm\langle n,\alpha\rangle|\geq\varepsilon_{\ell_{j}}^{4\tau},\ 0<|n|\leq N_{\ell_{j}}.
\end{split}
\end{equation*}
The two inequalities above imply that, for any $Y_j\in\Lambda_{j},$
\begin{equation}\label{202207242}
\begin{split}		
|A_j^{-1}Y_{j}(\cdot+\alpha)A_j-Y_{j}(\cdot)|_{\tilde{h}_j}\geq
\varepsilon_{\ell_{j}}^{12\tau}|Y|_{\tilde{h}_j}.
\end{split}
\end{equation}	
The inequality above, together with \eqref{202207241}, implies that
$\Lambda_{N}\subset C_{{\tilde{h}_j}}^{nre}
(N_{\ell_j},\varepsilon_{\ell_{j}}^{12\tau}).$

Note $\varepsilon_{\ell_{j}}^{12\tau}\leq(2|A|)^{-30}$ and by \eqref{202207156},
\begin{equation*}
\begin{split}		
|\Tilde{F}_{j}|_{{\tilde{h}_j}}\leq3\varepsilon^{k-57\tau}_{\ell_{j}}
<\{\varepsilon_{\ell_{j}}^{12\tau}\}^{61/30},
\end{split}
\end{equation*}	
then by Lemma \ref{implicittheorem} with $\varepsilon_{\ell_{j}}^{12\tau}$ in place of $\eta$ we know that there exist $Y\in C^{nre}_{{\tilde{h}_j}}(\varepsilon_{\ell_{j}}^{12\tau}),$ $F^{*}_{j}\in C^{re}_{{\tilde{h}_j}}(\varepsilon_{\ell_{j}}^{12\tau})$  such that
\begin{equation}\label{202201185}
e^{Y_{j}(x+\alpha)}(A_{{j}}e^{\Tilde{F}_{j}(x)})e^{-Y_{j}(x)}=A_{{j}}
e^{F^{*}_{j}(x)},
\end{equation}
with the estimates
\begin{equation}\label{202107223}
\begin{split}
|Y_{j}|_{\tilde{h}_j}
\leq2\varepsilon_{\ell_{j}}^{-12\tau}|\Tilde{F}_j|_{\tilde{h}_j}
\leq6\varepsilon_{\ell_{j}}^{(k-57\tau)-12\tau}=
6\varepsilon_{\ell_{j}}^{k-69\tau},
\end{split}
\end{equation}
and
\begin{equation}\label{202207223}
|F^{*}_{j}-\mathbb{P}_{re}\Tilde{F}_{j}|_{\tilde{h}_j}
\leq2\varepsilon_{\ell_{j}}^{-74\tau}|\Tilde{F}_{j}|_{\tilde{h}_j}^{2}.
\end{equation}
Thus
\begin{equation}\label{202207223a}
|F^*_{j}|_{\tilde{h}_j} \leq 2 |\Tilde{F}_{j}|_{\tilde{h}_j}\leq 6\varepsilon^{k-57\tau}_{\ell_{j}}.
\end{equation}

Let $B_{j}=e^{Y_j}$, then \eqref{202107223} implies that
\begin{equation*}
|B_j-I|_{h_j} \leq 12\varepsilon_{\ell_j}^{k-69\tau},
\end{equation*}
which, together with the estimates about $B^{(j)}$ in \eqref{202201180}, yields
\begin{equation*}
\begin{split}
|B^{(j+1)}|_{h_{j+1}}&=|B^{(j)}B_{j}|_{h_{j+1}} \leq\varepsilon_{\ell_{j}}^{-18\tau}(1+12\varepsilon_{\ell_j}^{k-69\tau})\\
&\leq2\varepsilon_{\ell_{j}}^{-18\tau}<\varepsilon_{\ell_{j+1}}^{-18\tau}.
\end{split}
\end{equation*}
That is the estimates about $B^{(j+1)}$ in \eqref{202201180} and \eqref{202306020}
also hold.

Note $\Lambda_{N}\subset C_{{\tilde{h}_j}}^{nre}
(N_{\ell_j},\varepsilon_{\ell_{j}}^{12\tau}),$ then for $F^{*}_{j}\in C^{re}_{{\tilde{h}_j}}(\varepsilon_{\ell_{j}}^{12\tau})$ we have
\begin{equation*}
F^{*}_{j}=\langle F^{*}_{j}\rangle+\mathcal{R}_{N_{\ell_{j}}}F^{*}_{j}.
\end{equation*}
Let
\begin{equation}\label{2022072231}
e^{F^{*}_{j}}=:e^{\langle F^{*}_{j}\rangle}e^{f_{j+1}},
\end{equation}
which, together with the equality above, yields
\begin{equation*}
\begin{split}		
|f_{j+1}|_{\frac{3}{4}h_{j}} &\leq   2|\mc{R}_{N>N_{\ell_{j}}}F^{*}_{j}|_{\frac{3}{4}h_{j}}
=2\sum_{|n|>N_{\ell_{j}}}|\widehat{F^{*}_{j}}(n)
|e^{|n| \frac{3}{4}h_{j}}\\
&\leq 2e^{- N_{\ell_{j}}\frac{1}{4}h_{j}}\sum_{|n|>N_{\ell_{j}}}|\widehat{F^{*}_{j}}(n)
|e^{|n| h_{j}}\\
&\leq 2\ell_{j}^{-5\tau\ell_{j}^{s}} \times6\varepsilon^{k-57\tau}_{\ell_{j}}\leq\varepsilon^{2k}_{\ell_{j+1}},
\end{split}
\end{equation*}
that is the estimates about $f_{j+1}$ in \eqref{202201180} holds with $j+1$ in place of $j.$

Let $A_{j+1}:=A_{j} e^{\langle F^{*}_{j}\rangle}.$ Then by combing with \eqref{202207223a}, we have
\begin{equation*}
\|A_{j+1}-A_j\| \leq 20\varepsilon^{k-57\tau}_{\ell_{j}}\|A_j\|,
\end{equation*}
and, if we expand $A_{j} e^{\langle F^{*}_{j}\rangle}$ for second order of $\langle F^{*}_{j}\rangle,$ we have
\begin{equation}\label{anier1}
\|A_j+\langle F^{*}_{j}\rangle-A_{j+1}\| \leq 2\|A_j\||F^{*}_{j}|^2_{\tilde{h}_j}
\leq72\|A_j\|\varepsilon_{\ell_{j}}^{2k-114\tau}.
\end{equation}
Moreover, \eqref{202207223}, \eqref{estimate_F}, and \eqref{202207155} yield
\begin{equation}\label{twof}
\|\langle F^{*}_{j}\rangle-\langle\Tilde{F}_{j}\rangle\|
\leq2\varepsilon_{\ell_{j}}^{-74\tau}|\Tilde{F}_{j}|_{\tilde{h}_j}^{2}
\leq72\varepsilon_{\ell_{j}}^{2k-188\tau},
\end{equation}
and
\begin{equation}\label{f_v}
\|\langle\tilde{F}_j(x)\rangle-\langle Ad(B^{(j)}) V_j(x) W\rangle\|\leq 2|f_j|_{\tilde{h}_j}\leq 2\varepsilon_{\ell_j}^{2k},
\end{equation}
respectively. Combing \eqref{anier1}-\eqref{f_v} with \eqref{202207223a}, we arrive at
\begin{equation*}
\|A_{j+1}-A_j-\langle Ad(B^{(j)}) V_j(x) W\rangle\|
\leq80 \varepsilon_{\ell_j}^{2k-188\tau},
\end{equation*}
which implies \eqref{3j1} and \eqref{3j2}.

\textbf{Case two:} $A_{j}\in \mathcal{RS}(N_{\ell_{j}},\varepsilon^{4\tau}_{\ell_{j}})$,  we only need to consider the case in which $A_j$ is elliptic with eigenvalues $\left\{e^{i \rho_j}, e^{-i \rho_j}\right\}$ for $\rho \in \mathbb{R} \backslash\{0\}$, because if $\rho_j \in i \mathbb{R}$, the non-resonant condition is always fulfilled due to the Diophantine condition on $\alpha$ and then it indeed belongs to the non-resonant case. For simplicity, set $0< \rho_j<2\pi$.
\begin{claim}
$n^{*}_{j}$ is the unique resonant site with
$ 0<|n^{*}_{j}| \leqslant N_{\ell_{j}}$.
\end{claim}
\begin{proof}
Indeed, if there exists $n_{j}^{\prime} \neq n^{*}_{j}$ satisfying $|2 \rho_{j}-\langle n_{j}^{\prime}, \alpha\rangle|<\varepsilon_{\ell_j}^{4\tau}$, then by the Diophantine condition of $\alpha$, we have
$$
\gamma|n_{j}^{\prime} - n^{*}_{j}|^{-\tau} \leq|\langle n_{j}^{\prime} - n^{*}_{j}, \alpha\rangle| <2 \varepsilon_{\ell_{j}}^{4\tau},
$$
which implies that $$|n_{j}^{\prime}|> \tilde{N}_{\ell_j}-N_{\ell_j} \gg  N_{\ell_{j}},$$
since $4/5<s<1,$ where
\begin{equation}\label{202305200}
\tilde{N}_{\ell_j}=(5^{-1}\gamma)^{\tau^{-1}} \varepsilon_{\ell_{j}}^{-4}.
\end{equation}
\end{proof}

The Diophantine condition of $\alpha$ and \eqref{202305200} yield
\begin{equation}\label{dc}
|\langle n, \alpha\rangle| \geq 5 \varepsilon_{\ell_j}^{4\tau} , \quad n \in \mathbb{Z}^d, 0<|n| \leq \tilde{N}_{\ell_j}.
\end{equation}
Then, the fact $|2 \rho_{j}-\langle n^{*}_{j}, \alpha\rangle|<\varepsilon^{4\tau}_{\ell_{j}}$ and the inequality above yield
\begin{equation*}
2|\rho_{j}| \geq|\langle n^{*}_{j}, \alpha\rangle|- \varepsilon^{4\tau}_{\ell_{j}} \geq 4 \varepsilon^{4\tau}_{\ell_{j}}.
\end{equation*}
which, together with Lemma 8.1 of \cite{HouY12} enable us to find $P_j \in S U(1,1)$ with the estimate
\begin{equation}\label{trans_P}
|P_j|^2 \leqslant 2\big(1+\frac{2}{|\rho_j|}\big) \leq3 \varepsilon^{-4\tau}_{\ell_j},
\end{equation}
such that
$$ A^{\prime}_j \triangleq P_j^{-1} A_j P_j=\left(\begin{array}{cc}e^{i \rho_j} & 0 \\ 0 & e^{-i \rho_j}\end{array}\right),  $$ where $0<\rho_j <2\pi$.
Set $g_{j} =P_{j}\Tilde{F}_{j}P^{-1}_{j}$, by \eqref{estimate_F} and \eqref{trans_P},
\begin{equation}\label{new_epsilon}
\left|g_{j}\right|_{\tilde{h}_j} \leq\|P_{j}\|^2|\Tilde{F}_{j}|_{\tilde{h}_j}
\leq9\varepsilon^{k-61\tau}_{\ell_j}.
\end{equation}
Denote by ${C}_{{\tilde{h}_j}}^{n r e}$ the set of all $Y=\left(\begin{array}{cc}i a(\theta) & y(\theta) \\ y(\theta) & -i a(\theta)\end{array}\right)$ in ${C}_{{\tilde{h}_j}}(\mathbb{T}^d, s u(1,1))$ with
$$
a(\theta)=\sum_{n \in \mathbb{Z}^d, 0<|n|<\tilde{N}_{\ell_j}} \widehat{a}(n) e^{i\langle n, \theta\rangle}, y(\theta)=\sum_{n \in \mathbb{Z}^d, 0 \leq|n|<\tilde{N}_{\ell_j}, n \neq n^*_{j}} \widehat{y}(n) e^{i\langle n, \theta\rangle},
$$
and by ${C}_{{\tilde{h}_j}}^{r e}$ the set of all $Z \in {C}_{{\tilde{h}_j}}(\mathbb{T}^d, s u(1,1))$ of the form
$$
Z(\theta)=\left(\begin{array}{cc}
i a & 0 \\
0 & -i a
\end{array}\right)+\left(\begin{array}{cc}
0 & \widehat{z}(n^*_j)e^{i\langle n^*_j, \theta\rangle} \\
\overline{\widehat{z}}(n^*_j) e^{-i\langle n^*_j, \theta\rangle} & 0
\end{array}\right)+\sum_{n \in \mathbb{Z}^d,|n| \geq \tilde{N}_{\ell_j}} \widehat{Z}(n) e^{i\langle n, \theta\rangle}
$$
where $a \in \mathbb{R}$, $\widehat{z}(n^*_j)\in\mathbb{C} $ and $\widehat{Z}(n)\in su(1,1)$ for $|n| \geq \tilde{N}_{\ell_j}.$

Applying Lemma \ref{implicittheorem}  to $(\alpha, A^{\prime}_{j} e^{g_j})$ to remove all the non-resonant terms of $g_j$ , then there exists $Y_{j} \in {C}_{{\tilde{h}_j}}$ and $g^{*}_{j} \in {C}_{{\tilde{h}_j}}^{r e}(\varepsilon_{\ell_j}^{4\tau})$ such that
$$
e^{Y_{j}(\theta+\alpha)}(A^{\prime}_{j} e^{g_{j}(\theta)}) e^{-Y_j(\theta)}=A^{\prime}_{j} e^{g^{*}_{j}(\theta)},
$$
with
\begin{equation}\label{two_g*}
|Y_j|_{{\tilde{h}_j}} \leq 2\varepsilon_{\ell_j}^{-4\tau} |g_j|_{{\tilde{h}_j}},\
|g^{*}_{j}-g_{j}^{r e}|_{{\tilde{h}_j}} \leq 2 \varepsilon_{\ell_j}^{-25\tau}|g_j|_{{\tilde{h}_j}}^2 .
\end{equation}
Denote $\breve{B}_{j} \triangleq e^{Y_j},$ which together with \eqref{two_g*} and \eqref{new_epsilon}, yields
\begin{equation}\label{3.55}
|\breve{B}_{j}-I|_{h_j} \leq 18\varepsilon^{k-65\tau}_{\ell_j}, \quad
|g^{*}_{j}|_{{\tilde{h}_j}}\leq2|g_j|_{\tilde{h}_j}\leq 18\varepsilon_{\ell_j}^{k-61\tau}.
\end{equation}
Note that  $g^{*}_j(\theta) \in {C}_{{\tilde{h}_j}}^{(r e)},$ then $g^{*}_j(\theta)$ is of the form
\begin{equation}\label{form}
g^{*}_j(\theta)=\left(\begin{array}{cc}i a_{j} & b^*_{j} e^{i\langle  n^{*}_{j}, \theta\rangle} \\ \bar{b}^*_{j} e^{-i\langle n^{*}_{j}, \theta\rangle} & -i a_{j}\end{array}\right)+\mathcal{R}_{{N}_{\ell_j}} g^{*}_{j},
\end{equation}
where $a_{j} \in \mathbb{R}, b^*_{j} \in \mathbb{C}$ with the estimates
\begin{equation}\label{est_c}
\begin{aligned}
&|a_{j}| \leq  |g^*_j|_{{\tilde{h}_j}} \leq 18\varepsilon_{\ell_j}^{k-61\tau}, \quad|b^*_j| \leq  |g^*_j|_{{\tilde{h}_j}} e^{-|n^{*}_j| {\tilde{h}_j}} \ll18\varepsilon_{\ell_j}^{k-61\tau}.
\end{aligned}
\end{equation}

Define the $4 \pi \mathbb{Z}^d$-periodic rotation $Q_{n^{*}_j}(\theta)$ :
$$
Q_{n^{*}_j}(\theta)=\left(\begin{array}{cc}
e^{-\frac{\langle n^*_j, \theta\rangle}{2} i} & 0 \\
0 & e^{\frac{\langle n^*_j, \theta\rangle}{2} i}
\end{array}\right),
$$
then we have
$$
Q_{n^{*}_j}(\theta+\alpha)\left(A^{\prime}_{j} e^{g^{*}_{j}(\theta)}\right) Q_{n^{*}_j}^{-1}(\theta)=\tilde{A}_{j} e^{g^{\prime}_j},
$$
with
$$
\tilde{A}_{j}=Q_{n^{*}_j}(\theta+\alpha) A^{\prime}_{j}Q_{n^{*}_j}^{-1}(\theta)=\left(\begin{array}{cc}
e^{i(\rho_j-\frac{\langle n^*_{j}, \alpha\rangle}{2})} & 0 \\
0 & e^{-i(\rho_j-\frac{\langle n^*_{j}, \alpha\rangle}{2})}
\end{array}\right),
$$
and
$$g^{\prime}_j= Q_{n^{*}_j} g^{*}_{j}(\theta) Q_{n^{*}_j}^{-1},
$$
which together with \eqref{form} yields that
$$
g^{\prime}_j=\tilde{C}_{j}  +\widetilde{g}_{j},
$$
where $\tilde{C}_{j} \triangleq\left(\begin{array}{cc}i a_{j} & b^{*}_{j} \\ \bar {b}^*_{j} & -i a_{j}\end{array}\right)$, $ \widetilde{g}_{j} \triangleq Q_{n^*_{j}} (\mathcal{R}_{{\tilde{N}_{\ell_j}}}g^{*}_{j}) Q_{n^*_{j}}^{-1} $.

Let $ B_{j}=Q_{n^*_{j}}  e^{Y_j}  P_{j} $,
 one can also show that
\begin{equation}
B_{j}(\theta+\alpha)(A_{j} e^{\Tilde{F}_{j}(\theta)}) B_{j}^{-1}(\theta)=\tilde{A}_j e^{\tilde{C}_{j}+\widetilde{g}_{j}}.
\end{equation}
Denote
\begin{equation}\label{3.53}
 e^{f_{j+1}(\theta)}=e^{-\tilde{C}_{j}}e^{\tilde{C}_{j}+\widetilde{g}_{j}},
\end{equation}
then by \eqref{202305200} and \eqref{3.55} we have
\begin{equation*}
\begin{aligned}
|f_{j+1}|_{\frac{3}{4} h_j}
& \leq|Q_{n^{*}_{j}}|_{\frac{3}{4} h_j}^2 \sum_{|n| \geq \tilde{N}_{\ell_j}}|g^{*}_{j}(n)| e^{\frac{3}{4}|n| h_j} \\
& \leq e^{\frac{3}{4}  N_{\ell_{j}} h_j} e^{-\frac{h_j}{4}\tilde{N}_{\ell_j}}|g^{*}_{j}|_{\widetilde{h}_j} \\
&=\exp\{-\frac{1}{4}h_j( \tilde{N}_{\ell_j}-3N_{\ell_{j}})\}|g^{*}_{j}|_{\widetilde{h}_j} \\
& \leq \exp\{-\frac{1}{8}h_j \tilde{N}_{\ell_j}\}
18\varepsilon_{\ell_j}^{k-61\tau}\ll \varepsilon_{\ell_{j+1}}^{2k},
\end{aligned}
\end{equation*}
that is \eqref{202201180} also holds in this case.

Let $ A_{j+1}=\tilde{A}_{j} e^{\tilde{C}_{j}}=e^{C_{j+1}}$, then recall the fact that $|2 \rho_{j}-\langle n^{*}_{j}, \alpha\rangle|<\varepsilon^{4\tau}_{\ell_{j}}$  and \eqref{3.55}, we have
\begin{equation}\label{6jj}
|C_{j+1}|\leqslant 2(|\rho_{j}-\frac{\langle n^*_{j}, \alpha\rangle}{2}|+|\widetilde{g}|_{h_j}) \leqslant \frac{5}{4}\varepsilon^{4\tau}_{\ell_{j}}.
\end{equation}
If denote $spec(C_{j+1})={\pm \mu_{j+1}}$ and ${C}_{j+1} \triangleq\left(\begin{array}{cc}i c_{j+1} & d_{j+1} \\ \bar {d}_{j+1} & -i c_{j+1}\end{array}\right)$, then $$\left|\mu_{j+1}\right| \leq \sqrt{\left.|| c_{j+1}\right|^2-|d_{j+1}|^2 \mid} \leq \frac{5\sqrt{2}}{4}\varepsilon^{4\tau}_{\ell_{j}},$$which implies that
$$|\rho_{j+1}|\triangleq |rot(\alpha,A_{j+1} )|=|\Im \mu_{j+1}| \leq 2\varepsilon^{4\tau}_{\ell_{j}}.$$
Moreover, \eqref{202305200} and \eqref{dc} imply that $\tilde{N}_{\ell_j}=5^{-\frac{1}{\tau}}
\gamma^{\frac{1}{\tau}} \varepsilon_{\ell_{j}}^{-4}>2\ell_{j+2},$ since $s<1,$ which, together with the inequality above, yields
\begin{equation*}
|2\rho_{j+1}-\langle n, \alpha\rangle| \geq
5\varepsilon_{\ell_{j}}^{4\tau}- 2|\rho_{j+1}|>\varepsilon^{4\tau}_{\ell_{j}},
\forall n \in \mathbb{Z}^d, 0<|n| \leq {N}_{\ell_{j+1}}.
\end{equation*}
Thus
\begin{equation*}
A_{j+1} \in \mathcal{N} \mathcal{R}\big(N_{\ell_{j+1}}, \varepsilon^{4\tau}_{\ell_{j+1}}\big).
\end{equation*}

Denote $spec (\title{C}_j)=\{\pm\lambda_{j}\}$, by \eqref{est_c}, we have
\begin{equation}\label{est_lam}
\lambda_{j}^2=|a_j|^2-|b_j^{*}|^2
\leq  |g^*_j|_{{\tilde{h}_j}}.
\end{equation}
Recall $ A_{j+1}=\tilde{A}_{j} e^{\tilde{C}_{j}},$ then by \eqref{202305300} in Lemma \ref{sS} we get
\begin{equation*}
b_{j+1}=b_{j}^{*}\lambda_j^{-1}sinh(\lambda_j)e^{i(\rho_j-\frac{\langle n^*_{j}, \alpha\rangle}{2})},
\end{equation*}
which, together with \eqref{3.55},\eqref{est_c} and \eqref{est_lam}, the fact $|2 \rho_{j}-\langle n^{*}_{j}, \alpha\rangle|<\varepsilon^{4\tau}_{\ell_{j}}, g_{j} =P_{j}\Tilde{F}_{j}P^{-1}_{j},$ and the inequalities $\frac{sinh(x)}{x}\leq 1+x^2, e^{x}\leq 1+2x, \forall |x|<1,$  yields
\begin{equation}\label{anier21}
\begin{split}
|b_{j+1}-b^*_j|&=b^*_j\big(\lambda_j^{-1}sinh(\lambda_j)e^{i(\rho_j-\frac{\langle n^*_{j}, \alpha\rangle}{2})}-1\big)\\
&\leq
|g^*_j|^{2}_{{\tilde{h}_j}} e^{-|n^{*}_j| {\tilde{h}_j}}
\leq(18)^{2}|\Tilde{F}_{j}|_{\tilde{h}_j}^{2}
\|P_j\|^{4} e^{-|n^{*}_j|\tilde{h}_j}.
\end{split}
\end{equation}
Moreover, the second inequalities in \eqref{two_g*} and the fact
$g_{j} =P_{j}\Tilde{F}_{j}P^{-1}_{j}$ give
\begin{equation}\label{two_g}
| \hat{g^{*}_{j}}(n^{*}_j)- \hat{g}_{j}(n^{*}_j)|_{\tilde{h}_j}
\leq2\varepsilon_{\ell_{j}}^{-25\tau}|\Tilde{F}_{j}|_{\tilde{h}_j}^{2}
\|P_j\|^{4} e^{-|n^{*}_j|\tilde{h}_j},
\end{equation}
and
\begin{equation}\label{g_v}
| g_j(x)- P_jAd(B^{(j)}) V_j(x) WP_j^{-1}|_{\tilde{h}_j}\leq 2|f_j|_{\tilde{h}_j}\|P_j\|^{2},
\end{equation}
respectively.
Combing \eqref{202201180}, \eqref{trans_P}, \eqref{new_epsilon}, \eqref{3.55} with \eqref{anier21}-\eqref{g_v},
we have
 \begin{equation*}
\left\{\begin{array}{l}
\big|b_{j+1}-\widehat{T}_j\big(n^{*}_j\big)\big| \leq 243 \varepsilon_{\ell_j}^{2k-147\tau} e^{-|n^{*}_j| h_j}, \quad \mathcal{K} \cap \mathcal{Z}_j \neq \emptyset ; \\
\big|b_{j+1}\big| \leq 243 \varepsilon_{\ell_j}^{2k-147\tau} e^{-\frac{3}{4}|n^{*}_j| h_{j-1}}, \quad \mathcal{K} \cap \mathcal{Z}_j=\emptyset,
\end{array}\right.
\end{equation*}
which gives the estimates in \eqref{7j}.

Let us estimate $|B^{(j+1)}|_{h_{j+1}}$ for $B^{(j+1)}=B^{(j)}B_j$ with $B_j=Q_{n^{*}_j} \breve{B}_{j} P_j^{-1}$. Note that
\begin{equation*}
\begin{aligned}
|Q_{n^{*}_j}|_{h_{j+1}}& \leq \exp\{2^{-1} N_{\ell_j} h_{j+1}\}
=\exp\{10\ell_{j+1}\tau\ell_{j+1}^{-1}\ln\ell_{j+1}\}\\
&=\ell_{j+1}^{10\tau}=\varepsilon^{-10\tau}_{\ell_{j+1}}.
\end{aligned}
\end{equation*}
Moreover, $A_{j-1} \in \mathcal{N} \mathcal{R}\big(N_{\ell_{j-1}}, \varepsilon^{4\tau}_{\ell_{j-1}}\big),$ then the estimate in \eqref{202306020} with $j$ in place of $j+1$ gives
$|B^{(j)}|_{h_j} \leq 2\varepsilon_{\ell_{j-1}}^{-18\tau},$ which, together with \eqref{trans_P}, \eqref{3.55} and the inequality above, yields
\begin{equation*}
\begin{aligned}
|B^{(j+1)}|_{h_{j+1}}\leq 2 \varepsilon^{-18\tau}_{\ell_{j-1}}  \times 3\varepsilon^{-4\tau}_{\ell_j}\times \varepsilon_{\ell_{j+1}}^{-10\tau}
\leq \varepsilon_{\ell_{j+1}}^{-18\tau}.
\end{aligned}
\end{equation*}
Thus the estimate about $B^{(j+1)}$ in \eqref{202201180} also holds with $j+1$ in place of $j+1$ provided $s>4/5.$
\end{proof}

\section{Some useful estimate}
We now assume that the fibered
rotation number $\rho(E)$ of the considered cocycle is  $\frac{1}{2}\langle n_{J},\alpha  \rangle$ for arbitrary but fixed $n_{J} \in \mathcal{K}$ and thus $\{n_{J}\} = \mathcal{K} \cap \mathcal{Z}_J$.
Let us introduce some notations.
Define the integer $j_0$ as
\begin{equation}\label{j_0}
j_0 \triangleq \min \big\{1 \leq j \leq J \mid A_j \in \mathcal{R} \mathcal{S}\big(N_{\ell_j}, \varepsilon^{4\tau}_{\ell_j}\big)\big\}.
\end{equation}
Define $\tilde{n}^{(j)}$ inductively for all $j \geq 1$, as
$$
\tilde{n}^{(j+1)}=\left\{\begin{array}{l}
\tilde{n}^{(j)}, \quad  A_j \in \mathcal{N} \mathcal{R}\big( N_{\ell_{j}}, \varepsilon^{4\tau}_{\ell_j}\big) \\
\tilde{n}^{(j)}+{n}^{*}_j, \quad  A_j \in \mathcal{R S}\big( N_{\ell_{j}}, \varepsilon^{4\tau}_{\ell_j}\big)
\end{array}\right.
$$
with $\tilde{n}^{(1)}=0.$ Moreover, set
$$
W_j \triangleq A d\big(B^{(j)}\big).W=B^{(j)} W\big(B^{(j)}\big)^{-1}, \ \  \text {with} \ \ W_1=W.
$$
 For $j_0<j \leq J$, we write $W_j$ as
\begin{equation}\label{w_j}
W_j=\left(\begin{array}{cc}
i u_j & w_j e^{-i\left\langle \tilde{n}^{(j)}, \cdot\right\rangle} \\
\overline{w_j} e^{i\left\langle\tilde{n}^{(j)}, \cdot\right\rangle} & -i u_j
\end{array}\right),
\end{equation}
together with quantities defined as follows:
\begin{equation}\label{de4.1}
\begin{aligned}
& \xi_j \triangleq\left|\left\langle w_j\right\rangle\right|, \quad \mathcal{M}_j \triangleq\left|w_j\right|_{h_j}+\left|u_j\right|_{h_j}, \\
& \mathfrak{m}_j \triangleq \sup _{n \in \mathbb{Z}^d,|n| \geq\left|\widetilde{n}^{(j)}\right|} \frac{1}{2}\left(\left|\widehat{w}_j(n)\right|+\left|\widehat{u}_j(n)\right|\right).
\end{aligned}
\end{equation}

These quantities are important for the following proof, especially the estimate about $\xi_j, j \geq j_0.$
We want to estimate $\xi_j$ step by step.
However, in every step of KAM, $\xi_{j+1}$ is influenced by the quantities $\mathcal{M}_j, \mathfrak{m}_j$ besides $\xi_j$.
Thus, we need to estimate $\xi_j, \mathcal{M}_j, \mathfrak{m}_j$ for $j \geq j_0$ in every step.
We first estimate ${n}^{*}_j$ and $\widetilde{n}^{(j)}$.

\begin{lemma}\label{6.1}
For all $j \geq 1$, we have
\begin{equation}\label{65}
|\widetilde{n}^{(j)}| \leq 5\ell_{j}/2.
\end{equation}
For $j>1$ and $A_j \in \mathcal{R S}\big( N_{\ell_{j}}, \varepsilon^{4\tau}_{\ell_j}\big)$, we have
\begin{equation}\label{66}
|\widetilde{n}^{(j)}| \leq5\ell_{j-1}/2,
\end{equation}
\begin{equation}\label{67}
|n_j^*|> 2\ell_{j},
\end{equation}
\begin{equation}\label{68}
|\widetilde{n}^{(j+1)}|>\frac{99}{50} \ell_{j}.
\end{equation}
As a consequence, for $j>1$ with $A_j \in \mathcal{R S}\big( N_{\ell_{j}}, \varepsilon^{4\tau}_{\ell_j}\big)$,
\begin{equation}\label{69}
\begin{aligned}
e^{|\widetilde{n}^{(j)}| h_j} \leq \varepsilon^{-\frac{1}{10}\tau}_{\ell_j}, \ \ e^{-|{n}^{*}_j| h_j} \leq \varepsilon^{20\tau}_{\ell_j}, \ \
e^{-|\widetilde{n}^{(j+1)}| h_j} \leq \varepsilon^{19.8\tau}_{\ell_j} .
\end{aligned}
\end{equation}
\end{lemma}

\begin{proof}
 We prove \eqref{65} inductively. It is obviously true for $j=1$. Assume that \eqref{65} is true for $j$. In the case $A_j \in \mathcal{N R} \big(N_{\ell_j}, \varepsilon^{4\tau}_{\ell_j}\big),$  \eqref{65} holds trivially since $\widetilde{n}^{(j+1)}=\widetilde{n}^{(j)}$. If  $A_j \in \mathcal{R} \mathcal{S}\big(N_{\ell_j}, \varepsilon^{4\tau}_{\ell_j}\big)$, we have
\begin{equation*}
\begin{split}
|\widetilde{n}^{(j+1)}|
\leq|\widetilde{n}^{(j)}|+|n^{*}_j|\leq 5\ell_{j}/2+2\ell_{j+1}  \leq 5\ell_{j+1}/2.
\end{split}
\end{equation*}

Now, we verify the estimates in \eqref{66}-\eqref{69}, thus fix our attentions on
the case $A_j \in \mathcal{R S}\big( N_{\ell_{j}}, \varepsilon^{4\tau}_{\ell_j}\big)$ with $j>1.$

The estimate in \eqref{6j2} implies  $A_{j-1} \in \mathcal{N R }\big( N_{\ell_{j-1}}, \varepsilon^{4\tau}_{\ell_{j-1}}\big)$, thus $\widetilde{n}^{(j)}=\widetilde{n}^{(j-1)}$, which, together with \eqref{65}, yields \eqref{66}. We will prove \eqref{67} by contraction.
Assume $|n^{*}_j|\leq 2\ell_{j}$, then by \eqref{2j} in Lemma \ref{nonresonantcase}, we have,
$$
\begin{aligned}
& \big|\rho_j-\rho_{j-1}\big| \leq 15|A_j-A_{j-1}|^{1 / 2} \leq  \varepsilon^{k/2-27\tau}_{\ell_{j-1}}, \\
&\big|2 \rho_j-\big\langle n^{*}_j, \alpha\big\rangle\big|<\varepsilon^{4\tau}_{\ell_{j}},
\end{aligned}
$$
thus by trigonometric inequality, we have
$$
\big|2 \rho_{j-1}-\big\langle n^{*}_j, \alpha\big\rangle\big|\leq 2 \varepsilon^{k/2-27\tau}_{\ell_{j-1}}+\varepsilon^{4\tau}_{\ell_{j}}
<\varepsilon^{4\tau}_{\ell_{j-1}},
$$
which contradicts with $A_{j-1} \in \mathcal{N R }\big( N_{\ell_{j-1}}, \varepsilon^{4\tau}_{\ell_{j-1}}\big)$.
Thus \eqref{67} holds.
The estimates in \eqref{66} and \eqref{67} yield $|\widetilde{n}^{(j)}|<\frac{1}{100}|n^{*}_j|,$ which, together with the fact $\widetilde{n}^{(j+1)}=\widetilde{n}^{(j)}+n^{*}_j,$ implies that \eqref{68} holds.
The estimates in \eqref{202207121} and \eqref{66}-\eqref{68} yield the three estimates in \eqref{69}, we omit the details.
\end{proof}
Recalling the definition of $j_0$ in \eqref{j_0} it is obvious that $\widetilde{n}^{(j)}=0$ for $j \leq j_0$ and $B^{(j)}$ is close to the identity for $1 \leq j<j_0$ as $j_0>1$. Denote $\|A_0\|=\sigma >1$, it follows that
\begin{lemma}\label{6.2}
 For $1 \leq j \leq j_0$
$$
\|A_j\| \leq 2 \sigma, \quad|W_j|_{h_j} \leq 2 , \quad\|[A_j,\langle W_j\rangle]\| \geq 1/2.
$$
\end{lemma}
\begin{proof}
The definition of $j_{0}$ implies that $A_{j} \in \mathcal{NR}\big( N_{\ell_{j}}, \varepsilon^{4\tau}_{\ell_{j}}\big), 1 \leq j \leq j_0.$
The relation $W_{j+1}=A d\big(B_j\big)W_j,$ together with the first inequality in \eqref{2j} yields
\begin{equation}\label{w1}
\left|W_{j+1}-W_j\right|_{h_j} \leq 24  \varepsilon^{k-69\tau}_{\ell_j}.
\end{equation}
Thus the first two inequalities in this lemma are from $\|A_0\|=\sigma,$
$\|W_0\|\leq 1,$ \eqref{w1}
and the two inequalities in \eqref{2j}.
 Moreover,
direct computations yield $\|[A_0,\langle W_0\rangle]\| =\|A_0 \|=\sigma$, then the last inequality in this lemma is a direct consequence of \eqref{2j} and \eqref{w1}.
\end{proof}
\begin{remark}
    In fact, we can prove inductively $\|A_j\| \leq 2 \sigma-2^{-j-1}$ for all $j\geq 1$ by using \eqref{2j} and \eqref{6jj}.
\end{remark}
We give the estimates about $\mathcal{M}_j, \mathfrak{m}_j, \xi_j$ below.
\begin{lemma}\label{6.3} We have,
\begin{equation}\label{202306022}
\begin{aligned}
\mathcal{M}_{j_0+1} \leq 24 \varepsilon_{\ell_{j_0}}^{-4\tau}, \ \ \mathfrak{m}_{j_0+1} \leq 24\varepsilon_{\ell_{j_0}}^{16\tau}, \ \ \xi_{j_0+1} \geq (80)^{-1}\varepsilon_{\ell_{j_0}}^{4\tau} .
\end{aligned}
\end{equation}
\end{lemma}

\begin{proof}
Notice that $A_{j_0} \in \mathcal{R S}( N_{\ell_{j_0}}, \varepsilon^{4\tau}_{\ell_{j_0}})$, by Lemma \ref{nonresonantcase}, we can find $B_{j_0}=Q_{-n^{*}_{j_0}} \breve{B}_{j_0} P_{j_0}^{-1}$ with  the following  estimates, see \eqref{6j},
\begin{equation}\label{611}
|\breve{B}_{j_0}-I|_{h_{j_0}} \leq 18\varepsilon^{k-65\tau}_{\ell_{j_0}}, \quad\|P_{j_0}\|^2 \leq 3\varepsilon^{-4\tau}_{\ell_{j_0}},\quad0<\rho_{j_0}<2\pi .
\end{equation}

Follow the notion in \eqref{w_j} and notice that $\widetilde{n}^{(j_0)}=0,$ we denote
\begin{equation*}
A d(\breve{B}_{j_0} P_{j_0}^{-1}) \cdot W_{j_0}
=\left(\begin{array}{cc}
i u_{j_0+1} & w_{j_0+1}  \\
w_{j_0+1} & -i u_{j_0+1}
\end{array}\right),
\end{equation*}
and
\begin{equation*}
W_{j_0+1}=A d(B_{j_0}W_{j_0})=
A d\left(Q_{-n^{*}_{j_0}}\right) \cdot\left(\begin{array}{cc}
i u_{j_0+1} & w_{j_0+1}  \\
w_{j_0+1} & -i u_{j_0+1}
\end{array}\right) .
\end{equation*}
Moreover, by telescoping, we obtain
\begin{equation}\label{detaiargu}
\begin{aligned}
|Ad(\breve{B}_{j_0} P_{j_0}^{-1})& \cdot W_{j_0}- A d( P_{j_0}^{-1}) \cdot W_{j_0}|_{h_{j_0}}\\
&\leq2|{\breve{B}_{j_0}}^{-1}|_{h_{j_0}}|\breve{B}_{j_0}-I|_{h_{j_0}}| Ad( P_{j_0}^{-1}) \cdot W_{j_0}|_{h_{j_0}}\\
&\leq2|{\breve{B}_{j_0}}^{-1}|_{h_{j_0}}
|\breve{B}_{j_0}-I|_{h_{j_0}}\|P_{j_0}\|^2|W_{j_0}|_{h_{j_0}},
\end{aligned}
\end{equation}
which, together with the fact $|W_{j_0}|_{h_{j_0}} \leq 2$  in Lemma \ref{6.2}
and \eqref{611}, yields
\begin{equation}\label{612}
\mathcal{M}_{j_0+1} \leq 4\|P_{j_0}\|^2(1+4|\breve{B}_{j_0}-I|_{h_{j_0}}) \leq 8\|P_{j_0}\|^2 \leq 24 \varepsilon_{\ell_{j_0}}^{-4\tau}.
\end{equation}
By Lemma \ref{6.2}, $\left|\left[A_{j_0},\left\langle W_{j_0}\right\rangle\right]\right| \geq \frac{1}{2}$. Then, using \eqref{3.31}, \eqref{611} and Lemma \ref{3.5}, we get
$$
\xi_{j_0+1} \geq \frac{1}{4 \rho_{j_0}\|P_{j_0}\|^2}-4|\breve{B}_{j_0}-I|_{h_{j_0}}
\|P_{j_0}\|^2|W_{j_0}|_{h_{j_0}} \geq (80)^{-1} \varepsilon_{\ell_{j_0}}^{4\tau} .
$$

We divide the estimate $\mathfrak{m}_{j_0+1}$ into two different cases. When $j_0>1$, by \eqref{69} and \eqref{612},
\begin{equation}\label{m}
\mathfrak{m}_{j_0+1} \leq \mathcal{M}_{j_0+1} e^{-|n^{*}_{j_0}| h_{j_0}} \leq
 24\varepsilon_{\ell_{j_0}}^{-4\tau}\varepsilon_{\ell_{j_0}}^{20\tau}  \leq 24\varepsilon_{\ell_{j_0}}^{16\tau}.
\end{equation}
\begin{remark}
Recall the definition in \eqref{de4.1}, $\mathcal{M}_{j_0+1} =\left|w_{j_0+1} \right|_{h_{j_0+1} }+\left|u_{j_0+1} \right|_{h_{j_0+1}}$, but in \eqref{detaiargu} and \eqref{612},  we actually bound $\mathcal{M}_{j_0+1}$ by $\left|w_{j_0+1} \right|_{h_{j_0} }+\left|u_{j_0+1} \right|_{h_{j_0} }$, yielding the first inequality in \eqref{m} to hold.
\end{remark}
 When $j_0=1$, recall the definition that $\tilde{n}^{(1)}=0$ and  \eqref{611}, similar to the argument in \eqref{detaiargu}, we have
$$
\mathfrak{m}_2 \leq|W_{+}-\langle W_{+}\rangle|_{h_1} \leq 4|\breve{B}_1-I|_{h_1}  \|P_1\|^2 \leq 216
\varepsilon_{\ell_1}^{k-65\tau}
\varepsilon_{\ell_{1}}^{-4\tau}\leq24\varepsilon_{\ell_{1}}^{16\tau} ,
$$
since $k>190\tau,$ where $W_{+} \triangleq A d\left(\breve{B}_1 P_1^{-1}\right) \cdot W=\left(\begin{array}{cc}i u_2 & w_2 \\ \overline{w_2} & -i u_2\end{array}\right)$ .
\end{proof}

\begin{lemma}\label{6.4}
Let $j_0+1 \leq j \leq J$. If $A_j \in \mathcal{N R} \big(N_{\ell_j}, \varepsilon^{4\tau}_{\ell_j}\big)$, then
\begin{equation}\label{613}
\begin{array}{ll}
\mathcal{M}_{j+1} \leq \mathcal{M}_j(1+24\varepsilon^{k-69.1\tau}_{\ell_j}), & \mathfrak{m}_{j+1} \leq \mathfrak{m}_j+24\mathcal{M}_j \varepsilon^{k-69.1\tau}_{\ell_j}, \\
\xi_{j+1} \geq \xi_j-24\mathcal{M}_j \varepsilon^{k-69.1\tau}_{\ell_j}.
\end{array}
\end{equation}
And if $A_j \in \mathcal{ RS} \big(N_{\ell_j}, \varepsilon^{4\tau}_{\ell_j}\big)$, then
\begin{equation}\label{614}
\begin{aligned}
& \mathcal{M}_{j+1} \leq 12 \mathcal{M}_j \varepsilon^{-4.1\tau}_{\ell_j}, \quad \mathfrak{m}_{j+1} \leq 12\mathcal{M}_j \varepsilon^{15.7\tau}_{\ell_j}, \\
& \xi_{j+1} \geq \xi_j-3 \mathfrak{m}_j-216\mathcal{M}_j \varepsilon^{k-69.1\tau}_{\ell_j}.
\end{aligned}
\end{equation}
\end{lemma}

\begin{proof}
If $A_{j} \in \mathcal{NR }\big( N_{\ell_{j}}, \varepsilon^{4\tau}_{\ell_{j}}\big)$,
by the first inequality in \eqref{2j} and \eqref{65}, we have the following  estimates
$$
\big|{B}_j-I\big|_{h_j} \leq 6\varepsilon^{k-69\tau}_{\ell_j}, \quad e^{|\widetilde{n}^{(j)}| h_j} \leq \varepsilon^{-\frac{1}{10}\tau}_{\ell_j} .
$$
It follows that
$$
\begin{aligned}
|W_{j+1}-W_j|_{h_j} \leq 4|B_j-I|_{h_j}|W_j|_{h_j}\leq 4|B_j-I|_{h_j}e^{|\widetilde{n}^{(j)}| h_j}|\mathcal{M}_j|_{h_j} \leq 24\varepsilon^{k-69.1\tau}_{\ell_j} \mathcal{M}_j,
\end{aligned}
$$
recall the definition in \eqref{de4.1}, thus the estimates in \eqref{613} hold.

If $A_{j} \in \mathcal{R S}\big( N_{\ell_{j}}, \varepsilon^{4\tau}_{\ell_{j}}\big),$ by \eqref{6j} in Lemma \ref{nonresonantcase}, we can find $B_{j}=Q_{-n^{*}_{j}} \breve{B}_{j} P_{j}^{-1}$ with certain estimates given in \eqref{611}.
Same as the proof of \eqref{detaiargu} in  Lemma \ref{6.3}, combining the estimates in \eqref{611} with \eqref{69} gives
$$
\begin{aligned}
\mathcal{M}_{j+1} \leq 2\|P_j\|^2(1+4| \breve{B}_j-I|_{h_j})|W_j|_{h_j} \leq 12\varepsilon^{-4.1\tau}_{\ell_{j}}\mathcal{M}_j .
\end{aligned}
$$
 By  Lemmas \ref{3.4} and \ref{3.5}, together with \eqref{69}, we have
$$
\begin{aligned}
\xi_{j+1} & \geq \xi_j-3 \mathfrak{m}_j-4 \mathcal{M}_j| \breve{B}_j -I|_{h_j}\times\|P_j\|^2 e^{|\widetilde{n}^{(j)}| h_j} \\
& \geq \xi_j-3 \mathfrak{m}_j-4 \mathcal{M}_j \times18\varepsilon^{k-65\tau}_{\ell_{j}} \times  3\varepsilon^{-4\tau}_{\ell_{j}}
\times\varepsilon^{-\frac{1}{10}\tau}_{\ell_{j}} \\
& \geq \xi_j-3 \mathfrak{m}_j-216\mathcal{M}_j \varepsilon^{k-69.1\tau}_{\ell_{j}}, \\
\mathfrak{m}_{j+1} & \leq \mathcal{M}_{j+1} e^{-|\widetilde{n}^{(j+1)}| h_j} \leq 12\mathcal{M}_j \varepsilon^{-4.1\tau}_{\ell_{j}} \varepsilon^{19.8\tau}_{\ell_{j}}=12\mathcal{M}_j \varepsilon^{15.7\tau}_{\ell_{j}}.
\end{aligned}
$$
\end{proof}

\begin{lemma}\label{6.5+}
For all $j_0+1\leq j\leq J,$ we have
\begin{equation}\label{6161}
\mathcal{M}_j \leq \varepsilon^{-6\tau}_{\ell_{j}}.
\end{equation}
\end{lemma}
\begin{proof}
We prove \eqref{6161} inductively. First,
by  Lemma \ref{6.3}, we have
\begin{equation*}
\mathcal{M}_{j_0+1} \leq 24\varepsilon^{-4\tau}_{\ell_{j_0}}<\varepsilon^{-6\tau}_{\ell_{{j_0+1}}},
\end{equation*}
thus \eqref{6161} holds for  $j_{0}+1$.
Assume that it holds for some $p \geq j_0+1$, now we prove \eqref{6161} holds for $p+1$.

If $A_{p} \in \mathcal{NR}(N_{\ell_{p}}, \varepsilon^{4\tau}_{\ell_{p}}),$  by  \eqref{613} in Lemma \ref{6.4}, we have
$$
\mathcal{M}_{p+1} \leq 2 \varepsilon^{-{6\tau}}_{\ell_{p}} \leq \varepsilon^{-{6\tau}}_{\ell_{p+1}}.
$$

If $A_{p} \in \mathcal{RS}(N_{\ell_{p}}, \varepsilon^{4\tau}_{\ell_{p}})$, the first inequality in \eqref{614} and the induction hypothesis yield
\begin{equation*}
\mathcal{M}_{p+1} \leq 12\mathcal{M}_{p}\times\varepsilon^{-{4.1\tau}}_{\ell_{p}} \leq \varepsilon^{-{6\tau}}_{\ell_{p+1}},
\end{equation*}
where we use the fact $s>4/5.$ Thus, the estimates in \eqref{6161} hold for all $j$ with $j_{0}+1\leq j\leq J.$

\end{proof}

\begin{lemma}\label{6.5}
    For all $j_0+1\leq j\leq J$,
\begin{equation}\label{615}
\xi_j \geq 10 \mathfrak{m}_j+\frac{1}{2^j} \xi_{j_0+1}.
\end{equation}
\end{lemma}

\begin{proof}

We prove \eqref{615} inductively. By using  Lemma \ref{6.3}, direct computation shows that  \eqref{615} holds for $j=j_0+1$. Assume that it holds for some $j \geq j_0+1$, now we prove \eqref{615} holds for $j+1$.

 If $A_{j} \in \mathcal{NR }\big( N_{\ell_{j}}, \varepsilon^{4\tau}_{\ell_{j}}\big)$,  by \eqref{613} and \eqref{6161}, we have
\begin{equation}\label{xiestimate}
\begin{aligned}
\xi_{j+1} & \geq \xi_j-24\mathcal{M}_j  \varepsilon^{k-69.1\tau}_{\ell_{j}}\geq \xi_j-\varepsilon^{k-76\tau}_{\ell_{j}} , \\
\mathfrak{m}_{j+1} & \leq \mathfrak{m}_j+24\mathcal{M}_j \varepsilon^{k-69.1\tau}_{\ell_{j}} \leq \mathfrak{m}_j+\varepsilon^{k-76\tau}_{\ell_{j}} .
\end{aligned}
\end{equation}
The estimate about $\xi_{j_0+1}$ in \eqref{202306022}, together with the smallness of $\varepsilon_{\ell_{j}},$ yields
\begin{equation*}
\begin{aligned}
\varepsilon^{k-76\tau}_{\ell_{j}}\leq  80\varepsilon^{k-70\tau}_{\ell_{j}}\times(80)^{-1}\varepsilon^{4\tau}_{\ell_{j}}
\leq(11\times 2^{j+1})^{-1} \xi_{j_0+1},
\end{aligned}
\end{equation*}
which, together with the estimates in \eqref{xiestimate}, and the induction hypothesis \eqref{615}: $\xi_j \geq 10 \mathfrak{m}_j+2^{-j} \xi_{j_0+1},$ implies
\begin{equation*}
\begin{aligned}
\xi_{j+1} & \geq \xi_j-\varepsilon^{k-76\tau}_{\ell_{j}}
\geq 10 \mathfrak{m}_j+2^{-j} \xi_{j_0+1}-\varepsilon^{k-76\tau}_{\ell_{j}} \\
& \geq 10(\mathfrak{m}_{j+1}-\varepsilon^{k-76\tau}_{\ell_{j}})+2^{-j} \xi_{j_0+1}-\varepsilon^{k-76\tau}_{\ell_{j}} \\
&=10 \mathfrak{m}_{j+1}+2^{-(j+1)} \xi_{j_0+1}+(2^{-(j+1)} \xi_{j_0+1}-11 \varepsilon^{k-76\tau}_{\ell_{j}}) \\
& \geq 10 \mathfrak{m}_{j+1}+2^{-(j+1)} \xi_{j_0+1}.
\end{aligned}
\end{equation*}

 If $A_{j} \in \mathcal{RS }\big( N_{\ell_{j}}, \varepsilon^{4\tau}_{\ell_{j}}\big)$, we have, by \eqref{614} and \eqref{6161},
\begin{equation}\label{xi2}
\begin{aligned}
\xi_{j+1} & \geq \xi_j-3 \mathfrak{m}_j-12\mathcal{M}_j \varepsilon^{k-69.1\tau}_{\ell_{j}} \geq \xi_j-3 \mathfrak{m}_j-\varepsilon^{k-76\tau}_{\ell_{j}}, \\
\mathfrak{m}_{j+1} & \leq 12\mathcal{M}_j \varepsilon^{15.7\tau}_{\ell_{j}} \leq 12\varepsilon^{9.7\tau}_{\ell_{j}}\leq (100\times2^{j})^{-1}\xi_{j_0+1}.
\end{aligned}
\end{equation}
The induction hypotheses, $\mathfrak{m}_j \leq \frac{1}{10} \xi_j$ and $2^{-j} \xi_{j_0+1} \leq \xi_j$ and the estimates in \eqref{xi2}, yield
\begin{equation*}
\begin{aligned}
\xi_{j+1} & \geq \frac{7}{10} \xi_j-\varepsilon^{k-76\tau}_{\ell_{j}}
\geq \frac{7}{10\times 2^j} \xi_{j_0+1}-\varepsilon^{k-76\tau}_{\ell_{j}} \\
&\geq\frac{1}{5\times 2^j} \xi_{j_0+1}+\frac{1}{ 2^{j+1}} \xi_{j_0+1}-\varepsilon^{k-76\tau}_{\ell_{j}}\\
&=\frac{1}{5\times 2^j} \xi_{j_0+1}-\varepsilon^{k-76\tau}_{\ell_{j}}+\frac{1}{ 2^{j+1}} \xi_{j_0+1}\\
& \geq\frac{1}{10\times 2^j} \xi_{j_0+1}+\frac{1}{ 2^{j+1}} \xi_{j_0+1}                                                                                                \\
&
\geq 10 \mathfrak{m}_{j+1}+\frac{1}{2^{j+1}} \xi_{j_0+1}.
\end{aligned}
\end{equation*}
\end{proof}

\section{Proof of Theorem \ref{p1}}
Consider the cocycle $(\alpha,S_E^V(x))$ with $S_E^V(x)$ defined by \eqref{11}, assume that the fiber rotation number $\rho(E)=\frac{1}{2}\langle n_{J},\alpha  \rangle$ for arbitrary but fixed $n_{J} \in \mathcal{K}.$
The hypothesis $|n|>\ell_{*},\forall n\in\mathcal{K},$ in \eqref{202207162} and definition of $\ell_{*}$ in \eqref{202307131} yield
\begin{equation*}
\varepsilon_{\ell_{0}}^{k-26\tau}<(2\|A_{0}\|)^{-30},
\end{equation*}
which enables us to apply Lemma \ref{nonresonantcase}  to the cocycle $(\alpha,M{S_{E}^{V(x)}} M^{-1})$. Inductively, there exists $B^{(j)}\in C_{h_{j}}^{\omega}(2\TT^{d},SU(1,1))$
which conjugates the cocycle $(\alpha,M{S_{E}^{V(x)}} M^{-1})$
 to $(\alpha, A_{j}e^{f_{j}}\prod_{p\geq j}e^{(V_{p}W_{j})})$,
where $W_j=\operatorname{Ad}(B^{(j)})W$ and
the rotation number $\tilde{\rho}_j$ of the above system is
\begin{equation}\label{6.20}
\tilde{\rho}_j=\frac{1}{2}\big\langle n_J-\widetilde{n}^{(j)}, \alpha\big\rangle.
\end{equation}
Note that $\rho_j=\operatorname{rot}(\alpha, A_j)$, combing  \eqref{wanfanhou}, \eqref{2022071611}, \eqref{202207121} with   \eqref{202201180} and \eqref{202207155} , we have

\begin{equation}\label{621}
\begin{aligned}
\big|\rho_j-\tilde{\rho}_j\big| & \leq 10 \big( \sup _{\theta \in \mathbb{T}^d}\big|f_j+\big(\sum_{p=j}^{\infty} v_p\big) W_j\big|\big)^{1/2}\\
& \leq 10 \big(\varepsilon^{2k}_{\ell_{j}}+\big(\sum_{p=j}^{\infty} \frac{1}{n^{k}_{p}}\big) 2\varepsilon_{\ell_{j}}^{-36\tau}\big)^{1/2}\\
& \leq 10\big(\varepsilon^{2k}_{\ell_{j}}+{4\varepsilon_{\ell_{j}}^{k-36\tau}} \big)^{1/2}
\leq2^{-1}\varepsilon_{\ell_{j}}^{4\tau},
\end{aligned}
\end{equation}
where the last inequality is by the fact $k>48\tau.$ The inequality above and \eqref{6.20} yield
\begin{equation}\label{6.22}
\big|2 \rho_j-\big\langle n_J-\widetilde{n}^{(j)}, \alpha\big\rangle\big| \leq \varepsilon_{\ell_{j}}^{4\tau}.
\end{equation}

In the following, we will prove that for all $j\geq J + 1$, $A_j$  are uniformly
bounded below and have zero rotation numbers which implies that $A_{\infty}$ is
hyperbolic or parabolic.

\begin{lemma}\label{6.6.0}
For $J=1$, we have $A_{1} \in \mathcal{RS }\big( N_{\ell_{1}}, \varepsilon^{4\tau}_{\ell_{1}}\big)$.
For $J>1$, we have either
\begin{equation}\label{625}
A_{J-1} \in \mathcal{RS }\big( N_{\ell_{J-1}}, \varepsilon^{4\tau}_{\ell_{J-1}}\big), \quad n^{*}_{J-1}
=n_J-\widetilde{n}^{(J-1)};
\end{equation}
or
\begin{equation}\label{626}
A_{J} \in \mathcal{RS }\left( N_{\ell_{J}}, \varepsilon^{4\tau}_{\ell_{J}}\right), \quad n^{*}_J=n_J-\widetilde{n}^{(J)}.
\end{equation}

\end{lemma}

\begin{proof}
 In case that $J=1$, we have $\widetilde{\rho}_1=\frac{1}{2}\left\langle n_1, \alpha\right\rangle$. By \eqref{6.22}
\begin{equation*}
A_{1} \in \mathcal{RS }( N_{\ell_{1}}, \varepsilon^{4\tau}_{\ell_{1}}), \quad n^{*}_1=n_1 .
\end{equation*}

Now we consider the case $J>1.$ If $|n_J-\widetilde{n}^{(J-1)}| \leq N_{\ell_{J-1}},$ then by \eqref{6.22}
\begin{equation*}\label{6251}
A_{J-1} \in \mathcal{RS }\left( N_{\ell_{J-1}}, \varepsilon^{4\tau}_{\ell_{J-1}}\right) , \quad n^{*}_{J-1}=n_J-\widetilde{n}^{(J-1)} .
\end{equation*}
Otherwise, if $\left|n_J-\widetilde{n}^{(J-1)}\right|>N_{\ell_{J-1}}$, then combing \eqref{202207146}, \eqref{2022071611},  \eqref{202207121}, with \eqref{202207121}, \eqref{65} and \eqref{6.22}, we have
$$
\begin{aligned}
\big|2 \rho_{J-1}-\langle n, \alpha\rangle\big| \geq & \big|\big\langle n-\big(n_J-\widetilde{n}^{(J-1)}\big), \alpha\big\rangle\big| \\
& -\big|2 \rho_{J-1}-\big\langle n_J-\widetilde{n}^{(J-1)}, \omega\big\rangle\big| \\
\geq & \gamma |{2\ell_{J+1}}|^{-\tau}- \varepsilon_{\ell_{J-1}}^{4\tau} \\
> &\varepsilon^{4\tau}_{\ell_{J-1}} , \quad n \in \mathbb{Z}^d,|n| \leq N_{\ell_{J-1}},
\end{aligned}
$$
since $4/5<s<1,$ which implies that $A_{J-1} \in \mathcal{NR }\big( N_{\ell_{J-1}}, \varepsilon^{4\tau}_{\ell_{J-1}}\big)$. Thus, $\widetilde{n}^{(J)}=\widetilde{n}^{(J-1)}$. By Lemma \ref{6.1} and the fact ${n_{J}} \in \mathcal{K} \cap \mathcal{Z}_J,$ we have
$$
\big|n_J-\widetilde{n}^{(J)}\big|=\big|n_J-\widetilde{n}^{(J-1)}\big| \leq \ell_{J+1}+ 5\ell_{J-1}\leq 2\ell_{J+1} ,
$$
which, together with \eqref{202207121} and  \eqref{6.22}, yields
\begin{equation*}\label{6261}
A_{J} \in \mathcal{RS}\big( N_{\ell_{J}}, \varepsilon^{4\tau}_{\ell_{J}}\big), \quad n^{*}_J=n_J-\widetilde{n}^{(J)}.
\end{equation*}

\end{proof}

\begin{lemma}\label{6.6}
\begin{equation*}
A_{J+1} \in \mathcal{NR }\big( N_{\ell_{J+1}}, \varepsilon^{4\tau}_{\ell_{J+1}}\big), \ \|A_{J+1}\| \geq n^{-(k+5\tau)}_{J}, \ \tilde{\rho}_{J+1}=0.
\end{equation*}
\end{lemma}
\begin{proof}
According to Lemma~\ref{6.6.0}, we divide the proof  into two different cases:

\noindent(1)$A_{J-1} \in \mathcal{RS }\big( N_{\ell_{J-1}}, \varepsilon^{4\tau}_{\ell_{J-1}}\big)$ with $J>1$. By \eqref{625},
\begin{equation}\label{627}
n^{*}_{J-1}=n_J-\widetilde{n}^{(J-1)}, \quad \widetilde{n}^{(J)}=\widetilde{n}^{(J-1)}+n^{*}_{J-1}=n_J,
\end{equation}
which, together with \eqref{6.20}, yields
\begin{equation}\label{628}
\tilde{\rho}_J=\frac{1}{2}\big\langle n_J-\widetilde{n}^{(J)}, \alpha\big\rangle=\frac{1}{2}\big\langle n_J-n_J, \alpha\big\rangle=0.
\end{equation}
By \eqref{67}, \eqref{6j2} and \eqref{627}, we get
\begin{equation}\label{629}
\big|n^{*}_{J-1}\big| \geq \frac{999}{1000}\big|n_J\big|, \quad A_{J} \in \mathcal{NR }\big( N_{\ell_{j}}, \varepsilon^{4\tau}_{\ell_{j}}\big).
\end{equation}

Let us write
$$A_J=\left(\begin{array}{cc} a_J & b_J \\ \overline{b_J} & \overline{a_J}\end{array}\right), \quad A_J+\left\langle V_J W_J\right\rangle \triangleq\left(\begin{array}{cc} a_{J,+} & b_{J,+} \\ \overline{b_{J,+}} & \overline{a_{J+}} \end{array}\right),$$
where $|a_J|^2-|b_J|^2=1, a_J, b_J \in \mathbb{C}$ and  $W_J=A d\left(B^{(J)}\right) \cdot W$ is given by
$$
W_J=\left(\begin{array}{cc}
i u_J & w_J e^{-i\left\langle\widetilde{n}^{(J)}, \cdot\right\rangle} \\
\overline{w_J} e^{i\left\langle\widetilde{n}^{(J)}, \cdot\right\rangle} & -i u_J
\end{array}\right)=\left(\begin{array}{cc}
i u_J & w_J e^{-i\left\langle n_J, \cdot\right\rangle} \\
\bar{w}_J e^{i\left\langle  n_J, \cdot\right\rangle} & -i u_J
\end{array}\right) .
$$
Recalling \eqref{202207161}, we have the conclusion
$$
\mathcal{K} \cap \mathcal{Z}_J \neq \emptyset \Rightarrow \mathcal{K} \cap \mathcal{Z}_{J-1}=\emptyset,
$$
and thus, together with \eqref{7j} and \eqref{629},
\begin{equation}\label{202306220}
\begin{aligned}
|b_J| &\leq e^{-\frac{3}{4}|n^{*}_{J-1}| h_{J-2}} \leq e^{-\frac{3}{4} \times \frac{999}{1000}|n_J| h_{J-2}} \\
& \leq e^{-\frac{1}{2}|\ell_J|(10\tau \ln \ell_{J-2})\ell_{J-2}^{-1}}\leq n^{-2k}_{J},
\end{aligned}
\end{equation}
where we use the fact $\ell_J>\ell_{*} >(2k/5)^{s^{-1}}$ (see \eqref{202307131}).

In view of $V_J W_J=\frac{1}
{n^{k}_{J}}\left(e^{i\left\langle n_J,
\theta\right\rangle}+e^{-i\left\langle n_J, \theta\right\rangle}\right) W_J$, it is obvious that
\begin{equation}\label{b+}
b_{J,+}=b_J+\frac{1}{n^{k}_{J}}\left(\left\langle w_J\right\rangle+\widehat{w}_J\left(2 n_J\right)\right) .
\end{equation}
Then, by Lemma \ref{6.3}, Lemma \ref{6.5} and \eqref{202306220} (note that $J>J-1 \geq j_0$),
$$
\begin{aligned}
\left|A_J+\left\langle V_J W_J\right\rangle\right| & \geq\left|b_{J,+}\right| \geq\frac{1}{n^{k}_{J}}\left(\left|\left\langle w_J\right\rangle\right|-\left|\widehat{w}_J\left(2 n_J\right)\right|\right)-\left|b_J\right| \\
& \geq \frac{1}{n^{k}_{J}}\left(\xi_J-\mathfrak{m}_J\right)-n^{-2k}_{J}  \geq \frac{9}{10} \frac{1}{n^{k}_{J}} \xi_J-n^{-2k}_{J}
\\
& \geq \frac{1}{2^{J+8}} \frac{1}{n^{k}_{J}}\varepsilon_{\ell_{j_0}}^{4\tau}
 \geq\frac{2}{n^{k}_{J}}\varepsilon_{\ell_{J}}^{5\tau}
 = 2n^{-(k+5\tau)}_{J}.
\end{aligned}
$$
Note that $k>190\tau$, by \eqref{3j1} in Lemma \ref{nonresonantcase},
$$
\|A_{J+1}\| \geq |b_{J+1}| \geq 2n^{-(k+5\tau)}_{J}- 80\varepsilon^{2k-188\tau}_{\ell_J}\geq n^{-(k+5\tau)}_{J} .
$$

The conclusions Lemma~\ref{6.6.0} show that
$A_{J} \in \mathcal{NR }\big( N_{\ell_{J}}, \varepsilon^{4\tau}_{\ell_{J}}\big)$, which, together with \eqref{628}, \eqref{621} and Lemma \ref{rotationnumberpers}, yields that $\widetilde{\rho}_{J+1}=\widetilde{\rho}_J=0,$ and $\left|\rho_{J+1}\right|<\varepsilon^{k/2-18\tau}_{\ell_{J+1}}.$ Thus,
for any $0<|n|\leq N_{\ell_{J+1}},$
\begin{equation*}
\big|2\rho_{J+1}-\langle n, \alpha\rangle\big|\geq\gamma|N_{\ell_{J+1}}|^{-\tau}
-2\varepsilon^{k/2-18\tau}_{\ell_{J+1}}
>\varepsilon^{4\tau}_{\ell_{J+1}}.
\end{equation*}
It follows that $A_{J+1} \in  \mathcal{NR }\big( N_{\ell_{J+1}}, \varepsilon^{4\tau}_{\ell_{J+1}}\big)$.

(2)$A_{J-1} \in \mathcal{NR }\big( N_{\ell_{J-1}}, \varepsilon^{4\tau}_{\ell_{J-1}}\big)$.
By \eqref{626} in Lemma~\ref{6.6.0}, we have
\begin{equation}\label{630}
n^{*}_J=n_J-\widetilde{n}^{(J)}.
\end{equation}
Moreover, by \eqref{6j} there is $P_J \in S U(1,1)$
such that
$A d\left(P_J^{-1}\right) \cdot A_J=\left(\begin{array}{cc} e^{i \rho_J} & 0 \\ 0 & e^{-i
\rho_J}\end{array}\right)$
satisfying
\begin{equation}\label{631}
\|P_J\|^2 \leq 3\varepsilon^{-4\tau}_{\ell_J}, \quad (0 < \rho_J \leq 2\pi).
\end{equation}
Then, $\operatorname{Ad}(P_J^{-1} B^{(J)}) \cdot(V_J W)=V_J \widetilde{W}_J$, where
\begin{equation*}
\widetilde{W}_J \triangleq A d(P_J^{-1} B^{(J)}) \cdot W=A d\left(P_J^{-1}\right) \cdot W_J \triangleq\left(\begin{array}{cc}
i \widetilde{u}_J & \widetilde{w}_J e^{-i\left\langle\widetilde{n}^{(J)}, \cdot\right\rangle} \\
\overline{\widetilde{w}_J }e^{i\left\langle\widetilde{n}^{(J)}, \cdot\right\rangle} & -i \widetilde{u}_J
\end{array}\right).
\end{equation*}
Write $V_J \widetilde{W}_J=\left(\begin{array}{cc}\square & T_J \\ \overline{T_J} & \square\end{array}\right)$. In view of $V_J=\frac{1}{n^{k}_J}\left(e^{i\left\langle n_J, \theta\right\rangle}+e^{-i\left\langle n_J, \theta\right\rangle}\right)$, then
\begin{equation}\label{632}
\widehat{T}_J(n_J-\widetilde{n}^{(J)})
=\frac{1}{n^{k}_J}(\langle\widetilde{w}_J\rangle+\widehat{\widetilde{w}}_J(2 n_J)).
\end{equation}

In case that $J>j_0 \geq 1$, by \eqref{630}, \eqref{632}, Lemmas \ref{3.4}, \ref{3.5}, \ref{6.3} and \ref{6.5}
$$
\begin{aligned}
\big|\widehat{T}_J\big(n_J-\widetilde{n}^{(J)}\big)\big| & \geq \frac{1}{n^{k}_J}\big(\big|\big\langle\widetilde{w}_J\big\rangle\big|
-\big|\widehat{\widetilde{w}}_J\big(2 n_J\big)\big|\big) \\
& \geq \frac{1}{n^{k}_J}\Big\{\frac{\|P_J\|^2+1}{2}\big(\xi_J-3 \mathfrak{m}_J\big)-\|P_J\|^2 \mathfrak{m}_J\Big\} \\
& \geq \frac{1}{4} \frac{1}{n^{k}_J}\xi_J \geq \frac{1}{2^{J+9}} \frac{1}{n^{k}_J} \varepsilon^{4\tau}_{\ell_{j_0}}\geq n^{-(k+5\tau)}_{J}.
\end{aligned}
$$

As for $J=j_0$, by  Lemma \ref{3.4},\ref{6.2} and \eqref{631}, we have
$$
\begin{aligned}
\left|\left\langle\widetilde{w}_{j_0}\right\rangle\right| & \geq \frac{1}{4\left|\rho_{j_0}\right|\|P_{j_0}\|^2} \geq \frac{1}{100}\varepsilon^{4\tau}_{\ell_{j_0}}, \\
\left|\widetilde{w}_{j_0}\right|_{h_{j_0}} & \leq\|P_{j_0}\|^2\left|W_{j_0}\right|_{h_{j_0}} \leq 6\varepsilon^{-4\tau}_{\ell_{j_0}}.
\end{aligned}
$$
Note that $|n_{j_0}| h_{j_0} \geq 10\tau \ln{\ell_{j_0}}$, thus by \eqref{632}, we have
$$
\begin{aligned}
|\widehat{T}_{j_0}(n_{j_0}-\widetilde{n}^{(j_0)})| & \geq \frac{1}{n^{k}_{j_0}}(|\langle\widetilde{w}_{j_0}\rangle|-
|\widetilde{w}_{j_0}|_{h_{j_0}} e^{-2|n_{j_0}| h_{j_0}}) \\
& \geq \frac{1}{n^{k}_{j_0}}\big(100^{-1}\varepsilon^{4\tau}_{\ell_{j_0}}
-6\varepsilon^{-4\tau}_{\ell_{j_0}}e^{-2|n_{j_0}| h_{j_0}}\big)\\
&\geq 200^{-1}n^{-(k+4\tau)}_{j_0}.
\end{aligned}
$$
In any case, by \eqref{7j} in Lemma \ref{nonresonantcase} and the fact $k>151\tau$, we get
$$
\begin{aligned}
\|A_{J+1}\|\geq \big|b_{J+1}\big|
&\geq\big|\widehat{T}_J\big(n_J-\widetilde{n}^{(J)}\big)\big|
-243\varepsilon_{\ell_j}^{2k-147\tau}\\
&\geq 200^{-1}n^{-(k+4\tau)}_{J}-243\varepsilon^{2k-147\tau}_{\ell_{J}}
\geq n^{-(k+5\tau)}_{J}.
\end{aligned}
$$
Moreover, $\tilde{\rho}_{J+1}=0$ follows from $n^{*}_J=n_J-\widetilde{n}^{(J)},$
and $A_{J+1} \in \mathcal{NR }\big(N_{\ell_{J+1}}, \varepsilon^{4\tau}_{\ell_{J+1}}\big)$ since $A_{J} \in \mathcal{RS }\big(N_{\ell_{J}}, \varepsilon^{4\tau}_{\ell_{J}}\big).$
\end{proof}

We furthermore have the following conclusion.

\begin{lemma}\label{12UP}
For $j \geq J+2,  A_{j} \in \mathcal{NR }\big( N_{\ell_{j}}, \varepsilon^{4\tau}_{\ell_{j}}\big)$ and $\widetilde{\rho}_j=0.$ Moreover,
\begin{equation}\label{upperb}
 n^{-(k+5\tau)}_{J}\leq|b_j| \leq n^{-(k-62\tau)}_{J}.
\end{equation}

\end{lemma}
\begin{proof}
We prove inductively that for all $j \geq J+2$
$$
 A_{j} \in \mathcal{NR }( N_{\ell_{j}}, \varepsilon^{4\tau}_{\ell_{j}}), \quad|b_j| \geq(1+\frac{1}{2^j}) n^{-(k+5\tau)}_{J}, \quad \widetilde{\rho}_j=0 .
$$

Consider the case $j=J+2$ first. The fact $\mathcal{K} \cap \mathcal{Z}_J \neq \emptyset$ implies that $\mathcal{K} \cap \mathcal{Z}_{J+1}=\emptyset$. By  \eqref{3j2} in Lemma \ref{nonresonantcase}, we get
\begin{equation*}
|b_{J+2}| \geq|b_{J+1}|-80\varepsilon^{2k-188\tau}_{\ell_J} \geq(1+\frac{1}{2^{J+1}}) n^{-(k+5\tau)}_{J} \geq(1+\frac{1}{2^{J+2}}) n^{-(k+5\tau)}_{J}.
\end{equation*}
Note $A_{J+1} \in \mathcal{NR }\big( N_{\ell_{J+1}}, \varepsilon^{4\tau}_{\ell_{J+1}}\big),$ then combing \eqref{2j}, \eqref{621}   with  Lemma \ref{rotationnumberpers}, we have $\tilde{\rho}_{J+2}=0$ and $\big|\rho_{J+2}\big|<30\varepsilon^{k/2-18\tau}_{\ell_{J+2}}.$ Thus,
for any $0<|n|\leq N_{\ell_{J+2}}$,
\begin{equation}\label{non}
\big|2\rho_{J+2}-\langle n, \alpha\rangle\big|\geq\gamma|N_{\ell_{J+2}}|^{-\tau}
-60\varepsilon^{k/2-18\tau}_{\ell_{J+2}}
>\varepsilon^{4\tau}_{\ell_{J+2}},
\end{equation}
which implies that $A_{J+2} \in  \mathcal{NR }\big( N_{\ell_{J+2}}, \varepsilon^{4\tau}_{\ell_{J+2}}\big)$.

We inductively assume that the desired conclusion holds for $j \geq J+2$, and we now verify it for $j+1$. In fact, by \eqref{2j} in Lemma \ref{nonresonantcase} and  $|b_j|\leq|A_j|\leq 2\sigma$, we have
$$
\begin{aligned}
|b_{j+1}| & \geq|b_j|- 20|A_j|\varepsilon^{k-57\tau}_{\ell_j} \\
&\geq(1-40\sigma \varepsilon^{k-57\tau}_{\ell_j}) \times(1+2^{-j})   n^{-(k+5\tau)}_{J} \\
& \geq(1+2^{-(j+1)})  n^{-(k+5\tau)}_{J} .
\end{aligned}
$$
Note that $A_{j} \in \mathcal{NR }\big( N_{\ell_{j}}, \varepsilon^{4\tau}_{\ell_{j}}\big)$, combing  \eqref{2j}   with  Lemma \ref{rotationnumberpers},
  we have $\widetilde{\rho}_{j+1}=\tilde{\rho}_j=0$. Then \eqref{621} implies $|\rho_{j+1}|<30\varepsilon^{k/2-18\tau}_{\ell_{j+1}} $. Similar to \eqref{non}, we have $A_{j+1} \in \mathcal{NR }\big( N_{\ell_{j+1}}, \varepsilon^{4\tau}_{\ell_{j+1}}\big)$.

 For the upper bounds of $b_j$, similar to Lemma \ref{6.6}, we  also divide the proof into two cases:

(1)$A_{J-1} \in \mathcal{RS }\big( N_{\ell_{J-1}}, \varepsilon^{4\tau}_{\ell_{J-1}}\big)$. By \eqref{202306220} and \eqref{2j}, we have
\begin{equation}
\begin{aligned}
|b_j| &\leq |b_{J}|+\Sigma^{j-1}_{p=J} |b_{p+1}-b_{p}| \\
&\leq n^{-2k}_{J}+\Sigma^{j-1}_{p=J} |A_{p+1}-A_p|  \\
& \leq   n^{-2k}_{J}+\varepsilon^{k-58\tau}_{\ell_{J}}\leq n^{-k+62\tau}_{J}.
\end{aligned}
\end{equation}

(2)$A_{J-1} \in \mathcal{NR }\big( N_{\ell_{J-1}}, \varepsilon^{4\tau}_{\ell_{J-1}}\big)$. Recall the definition of $W$ and $V_j$, by \eqref{202201180} and \eqref{6j}, we have
\begin{equation*}
|P_{J-1} A d(B^{(J-1)}) \cdot(V_{J-1} W) P_{J-1}^{-1}|_{h_{J-1}} \leq 3 n_{J}^{-k+61\tau},
\end{equation*}
 then the inequality in right hand of \eqref{upperb} is a direct result of \eqref{7j} in Lemma \ref{nonresonantcase}.
\end{proof}

For $j \geq J+1,$ the conclusions of Lemmas~\ref{6.6},\ref{12UP} demonstrate that $A_{j} \in \mathcal{NR}( N_{\ell_{j}}, \varepsilon^{4\tau}_{\ell_{j}})$.  Recall  the inequalities \eqref{2j} and \eqref{202201180} in  Lemma~\ref{nonresonantcase},
 \begin{equation}\label{202307140}
\begin{aligned}
|B_j-I|_{h_j} \leq 6\varepsilon^{k-69\tau}_{\ell_j},
|B^{(J+1)}|_{h_{J+1}}\leq \varepsilon_{\ell_{J+1}}^{-18\tau}.
\end{aligned}
\end{equation}
Then, for any $j\in\NN$ with $j\geq J+2,$ we have
\begin{equation}
\begin{aligned}
|B^{(j)}|_{h_{j-1}}&= |B_{j-1}B_{j-2}\cdots B_{J+1}B^{(J+1)}|_{h_{j-1}}\\
&=(\prod_{p=J+1}^{j-1}|B_{p}|_{h_{p}})|B^{(J+1)}|_{h_{J+1}} \\
&\leq\prod_{p=J+1}^{j-1}(1+6\varepsilon^{k-69\tau}_{\ell_p})
\varepsilon_{\ell_{J+1}}^{-18\tau}\leq2\varepsilon_{\ell_{J+1}}^{-18\tau}.
\end{aligned}
\end{equation}
By Cauchy estimate, for $k_0 \in \mathbb{N}$ with $k_0 \leq k-69\tau, \forall j\geq J+1,$ we have
\begin{equation*}
\begin{aligned}
\|D^{k_{0}}(B^{(j+1)}-B^{(j)})\|_{C^0}
& \leq (e^{-1}k_{0})^{k_{0}}h_{j}^{-k_{0}}|B^{(j+1)}-B^{(j)}|_{h_{j}}\\
& \leq (e^{-1}k_{0})^{k_{0}}h_{j}^{-k_{0}}|B^{(j)}|_{h_{j}}
|B_{j}-I|_{h_{j}}\\
&\leq 12(10\tau ek_{0}^{-1}\ln\ell_{j})^{-k_{0}}\varepsilon_{\ell_{J+1}}^{-18\tau},
\end{aligned}
\end{equation*}
which, together with the inequality above, yields
\begin{equation}\label{cauchys}
\begin{aligned}
\|B^{(j+1)}-B^{(j)}\|_{k_{0}}
\leq (k_{0}^{-1}\ln\ell_{j})^{-k_{0}}\varepsilon_{\ell_{J+1}}^{-18\tau}.
\end{aligned}
\end{equation}

Set $B^{\infty}=\lim _{j \rightarrow \infty} B^{(j)}.$ We can deduce from \eqref{202307131} \eqref{2022071611} and \eqref{202207162} that $ n_J\geq\ell_{J}\geq e^{k},$ together with \eqref{202307140} and \eqref{cauchys}, yields
\begin{equation}\label{cauchy}
\begin{aligned}
\|B^{\infty}\|_{k_{0}}&\leq\|B^{(J+1)}\|_{k_{0}}+\sum_{j=J+1}^{\infty}
\|B^{(j+1)}-B^{(j)}\|_{k_{0}}\\
&\leq 2k_{0}^{k_{0}}\varepsilon_{\ell_{J+1}}^{-18\tau}\leq n_{J}^{k_0+36\tau},
\end{aligned}
\end{equation}
thus $B^{\infty}\in C^{k_0}(2\TT^{d},SU(1,1))$  conjugates the system $(\alpha,MS_E^V(x)M^{-1} )$ to $(\alpha, A_{\infty})$, where $A_{\infty}=\lim _{j \rightarrow \infty} A_j$.  As we have proved
 $\rho((\alpha, A_{\infty}))=0$ in Lemma~\ref{12UP},
 the eigenvalues of $A_{\infty}$ can not be $e^{i \theta}$ for some $\theta \in \RR$. Note that  $2 \rho(\alpha, S_E^V(x))-\langle n_J, \alpha\rangle \in \mathbb{Z}$ for $n_J \in \mathbb{Z}^d \backslash\{0\}$, and $(\alpha, A_0+F_0(\cdot))$ is not uniformly hyperbolic for $E$ in the spectrum. Assume that $A_{\infty}=\left(\begin{array}{cc}
* & b_{\infty} \\
\bar{b_{\infty}} & *
\end{array}\right)$ with $b_{\infty}=\lim _{j \rightarrow \infty} b_j$, then by Lemma \ref{RD}, there exists $R_{\phi}$  such that $R_{-\phi}M ^{-1}A_{\infty} MR_{\phi}=C$
where $C=\left(\begin{array}{cc}
1 & \zeta \\
0 & 1
\end{array}\right)$ with $\zeta= |b_{\infty}|$.
By taking $B=M^{-1}B^{\infty}MR_{\phi}$, we have
$$
B(\cdot+\alpha)^{-1}A(\cdot) B(\cdot)=\left(\begin{array}{cc}
1 & \zeta\\
0 & 1
\end{array}\right),
$$
 and \eqref{est12} can be derived directly from \eqref{upperb} and \eqref{cauchy}.

\section{Proof of Theorem \ref{gape}}
For any $  n_J\in\mathcal{K}$, we will estimate the size of
the spectral gap $I_{J}(V)=(E_{n_J}^{-},E_{n_J}^{+}).$ If  the energy of  Schr\"{o}dinger cocycle lies at right edge point of the gap,
by Theorem \ref{p1},  there exists
 $B\in C^{k_0}(2\TT^{d},SL(2,\RR))$ conjugates the system $(\alpha,S_{E_{n_J}^{+}}^V)$ to $(\alpha, C)$, that is
\begin{equation}\label{202307143}
B(\theta+\alpha)^{-1}S^{V}_{E_{n_J}^{+}}B(\theta) =C=\begin{pmatrix}
1&\zeta\\
0&1
\end{pmatrix},
\end{equation}
with the estimates :
\begin{equation}\label{202307141}
\|B\|_{k_0}\leq n_J^{k_0+36\tau},
\end{equation}
and
\begin{equation}\label{202307142}
n^{-(k+5\tau)}_{J}\leq|\zeta|\leq  n^{-(k-62\tau)}_{J}.
\end{equation}

In this section, we'll demonstrate how $|B|_{k_0}$ and $\zeta$ determine $|I_J(V)|$. To this end, we'll first present the $C^k$ variant of the Moser-P\"oschel argument, which was first established in \cite{CAI2021109035}.

\subsection{Moser-P\"{o}schel argument}

Assume that $\zeta\in(0,\frac{1}{2})$.
For any $\delta\in (0,1)$, by applying \eqref{202307143}, direct calculations yield
\begin{equation}\label{7mb}
B(\theta+\alpha)^{-1}S^{V}_{E_{n_J}^{+}-\delta}B(\theta) = C-\delta P(\theta),
\end{equation}
where
\begin{equation}
P(\theta) = \begin{pmatrix}
B_{11}(\theta)B_{12}(\theta)-\zeta B_{11}^{2}(\theta)&-\zeta B_{11}(\theta)B_{12}(\theta) +B_{12}^{2}(\theta)\\
-B_{11}^{2}(\theta)&-B_{11}(\theta)B_{12}(\theta)
\end{pmatrix}, \label{P}
\end{equation}
with estimate
\begin{equation}\label{7m6.2}
\|P(\theta)\|_{k_{0}} \leq (1+\zeta)\|B\|^{2}_{k_{0}},  \ k_{0}\leq k-90\tau.
\end{equation}

\begin{lemma}\emph{[Lemma 4.1 \cite{CAI2021109035}]}\label{kam}
	Suppose that $\alpha\in {\rm DC}_{d}(\gamma,\tau)$ and $P(\theta)\in C^{k_0}(2\TT^{d},SL(2,\RR))$ of form {\rm (\ref{P})}.
	Let $D_{\tau} = 8\sum_{m=1}^{\infty} (2\pi m)^{-(k_{0}-\hat{k}-3\tau-d+1)}$ with $\hat{k}\in
	\mathbb{Z}$ and $\hat{k}<k_{0}-3\tau-d$.
	If $0<\delta <  D_{\tau}^{-1}\gamma^{3}|B(\theta)|_{k_{0}}^{-2}$, then there exist $\widetilde{B}(\theta) \in C^{\hat{k}}(2\mathbb{T}^{d},{\rm SL}(2,\mathbb{R}))$  and $P_{1}(\theta)\in C^{\hat{k}}(\mathbb{T}^{d},{\rm gl}(2,\mathbb{R}))$ such that
	\begin{equation*}
	\widetilde{B}(\theta+\alpha)^{-1}(C-\delta P(\theta))\widetilde{B}(\theta) = e^{c_{0}-\delta c_{1} }+\delta^{2}P_{1}(\theta),
	\end{equation*}
	where $c_{0}=\begin{pmatrix}
	0&\zeta\\
	0&0
	\end{pmatrix}$ and
	\begin{equation} \label{b}
	c_{1}=\begin{pmatrix}
	[B_{11}B_{12}]-\frac{\zeta}{2}[B_{11}^{2}] &-\zeta[B_{11}B_{12}]+[B_{12}^{2}]\\
	-[B_{11}^{2}]&-[B_{11}B_{12}]+\frac{\zeta}{2}[B_{11}^{2}]
	\end{pmatrix}
	\end{equation}
	with estimates
	\begin{align*}
	&\|\widetilde{B}(\theta)-{\rm Id}\|_{\hat{k}} \leq D_{\tau} \gamma^{-3}\delta |B|_{k_{0}}^{2},\\
	&\|P_{1}(\theta)\|_{\hat{k}} \leq 53D_{\tau}^{2}\gamma^{-6}\|B\|_{k_{0}}^{4}  +\delta^{-1}\zeta^{2}\|B\|_{k_{0}}^{2}.
	\end{align*}
\end{lemma}

In the following, we  can obtain bounds on the gaps' lengths  based on the information that quantitative reducibility provides.

Since $\tilde{B}$ is homotopic to identity by construction, we have
$$
\rho(\alpha, C-\delta P(\cdot))=\rho\left(\alpha, e^{c_0-\delta c_1}+\delta^2 P_1(\cdot)\right).
$$
Following that, we will present a important quantity in our estimations of the length of the gaps. Specifically, we define a function $d(\delta):=\det(c_{0}-\delta c_{1}) + \frac{1}{4}\delta^{2}\zeta^{2}[B_{11}^{2}]^{2}$ for any $\delta \in(0,1)$.

A direct calculation can provide that
\begin{equation}
\begin{split}
d(\delta) & = -\delta [B_{11}^{2}]\zeta +\delta^{2}([B_{11}^{2}][B_{12}^{2}]-[B_{11}B_{12}]^{2})\\
& = \delta ([B_{11}^{2}][B_{12}^{2}]-[B_{11}B_{12}]^{2})\big(\delta - \frac{[B_{11}^{2}]\zeta }{[B_{11}^{2}][B_{12}^{2}]-[B_{11}B_{12}]^{2}} \big).
\end{split} \label{d}
\end{equation}
To  estimate $d(\delta)$ further, we recall the following fundamental lemma which was established in the $C^{\omega}$ case in Lemma 6.2 and Lemma 6.3 of \cite{zhaozhiyan} and it also holds in the $C^k$ case in Lemma 4.2 of \cite{CAI2021109035}.
\begin{lemma}[\cite{zhaozhiyan,CAI2021109035}] \label{poly}
For any $B\in C^{k_{0}}(2\mathbb{T}^{d},{\rm SL}(2,\mathbb{R}))$,  we have $$[B_{11}^2] \geq(2\|B\|_{k_0})^{-2}.$$
Furthermore, if $|B|_{k_{0}}\zeta^{\frac{\kappa}{2}}\leq \frac{1}{4}$ with $\kappa \in (0,\frac{1}{4})$,  the followings hold:
	\begin{align*}
	&0<\frac{[B_{11}^{2}]}{[B_{11}^{2}][B_{12}^{2}]-[B_{11}B_{12}]^{2}} \leq \frac{1}{2}\zeta^{-\kappa},\\
	&[B_{11}^{2}][B_{12}^{2}]-[B_{11}B_{12}]^{2} \geq 8\zeta^{2\kappa}.
	\end{align*}
\end{lemma}

\begin{proof}[Proof of Theorem \ref{gape}]

For any $n_J\in \mathbb{Z}^{d}\backslash\{0\}$, by \eqref{202307141} and \eqref{202307142}, we have
\begin{equation}\label{tiaojian}
\|B\|_{k_{0}}^{14}\zeta^{\frac{1}{10}}\leq  |n_J|^{14(k_{0}+36\tau)}  |n_J|^{\frac{-k+62\tau}{10}}\leq 10^{-11}D_{\tau}^{-4}\gamma^{12},
\end{equation}
the above inequality is possible since one can  fix $k_{0}=[\frac{k}{1000}]$ when $k>1000\tau$ is large enough.

Let $\delta_{1} = \zeta^{\frac{9}{10}}$.   Note that $|B|_{k_{0}}\geq 1$, by (\ref{tiaojian}), we have
\begin{equation*}
\delta_{1}D_{\tau} \gamma^{-3} \|B(\theta)\|_{k_{0}}^{2}\leq \zeta^{\frac{9}{40}} D_{\tau} \gamma^{-3}\|B(\theta)\|_{k_{0}}^{\frac{7}{2}} \leq 10^{-\frac{11}{4}}<1,
\end{equation*}
which deduces that $0<\delta_{1}<D_{\tau}^{-1}\gamma^{3}\|B(\theta)\|_{k_{0}}^{-2}$. Then by Lemma \ref{kam}, there exist $\widetilde{B}\in C^{\hat{k}}(2\mathbb{T}^{d},{\rm SL}(2,\mathbb{R}))$ and $P_{1}\in C^{\hat{k}}(\mathbb{T}^{d},{\rm gl}(2,\mathbb{R}))$ such that the cocycle $(\alpha, B-\delta_{1}P(\theta))$ is conjugated to $(\alpha, e^{c_{0}-\delta_{1}c_{1}}+\delta_{1}^{2}P_{1})$ by $\widetilde{B}$. Since $\widetilde{B}$ is homotopic to identity by construction, we have
$$
\rho(\alpha, C-\delta P(\cdot))=\rho\big(\alpha, e^{c_0-\delta c_1}+\delta^2 P_1(\cdot)\big).
$$
It is sufficient to demonstrate that
$\rho(\alpha, e^{c_{0}-\delta_{1}c_{1}}+\delta_{1}^{2}P_{1}(\theta))>0 $ and by monotonicity of rotation number, we have  $|I_J(V)|\leq \delta_{1}$.
 The equation \eqref{tiaojian} claims that $
\|B\|_{k_{0}}\zeta^{\frac{1}{20}} \leq 10^{-\frac{11}{2}}D_{\tau}^{-2} \|B\|^{-6}_{k_{0}} \gamma^{6} \leq \frac{1}{4}$.
Applying Lemma \ref{poly} to (\ref{d}),  we have
\begin{equation*}
d(\delta_{1}) \geq \zeta^{\frac{9}{10}}\times 8\zeta^{\frac{1}{5}}\times \frac{1}{2}\zeta^{\frac{9}{10}} =4\zeta^{2}.
\end{equation*}
Furthermore, by (\ref{tiaojian}) and $d(\delta)=\det(c_{0}-\delta c_{1}) + \frac{1}{4}\delta^{2}\zeta^{2}[B_{11}^{2}]^{2}$, it is straightforward to observe that
\begin{equation}
\begin{split}
\det (c_{0}-\delta_{1}c_{1})
\geq 4\zeta^{2}-\frac{1}{4}\delta_{1}^{2}\zeta^{2}[B_{11}^{2}]^{2}
\geq 4\zeta^{2}(1-\frac{1}{16}\zeta^{\frac{9}{5}}\|B\|_{k_{0}}^{4})\geq 3\zeta^{2}.
\end{split} \label{c1}
\end{equation}
Hence, by Lemma 8.1 of \cite{HouY12}, there exists $\mathcal{P}\in {\rm SL}(2,\mathbb{R})$ such that
\begin{equation*}
\mathcal{P}^{-1}e^{c_{0}-\delta_{1}c_{1}}\mathcal{P} = \exp \begin{pmatrix}
0&\sqrt{\det (c_{0}-\delta_{1}c_{1})}\\
-\sqrt{\det (c_{0}-\delta_{1}c_{1})}&0
\end{pmatrix}:=\Delta
\end{equation*}
with $
|\mathcal{P}|\leq 2 \left(\frac{|c_{0}-\delta_{1}c_{1}|}{\sqrt{\det (c_{0}-\delta_{1}c_{1})}}\right)^{\frac{1}{2}}$. Since $|c_{0}-\delta_{1}c_{1}|\leq \zeta+\delta_{1}(1+\zeta)\|B\|_{C^{0}}^{2}\leq \sqrt{3}\zeta^{\frac{9}{10}} \|B\|^{2}_{k_{0}}$, we have
\begin{equation*}
\begin{split}
\frac{|c_{0}-\delta_{1}c_{1}|}{\sqrt{\det (c_{0}-\delta_{1}c_{1})}}
\leq \frac{\sqrt{3}\zeta^{\frac{9}{10}}\|B\|^{2}_{k_{0}}}{\sqrt{3}\zeta} \leq \|B\|_{k_{0}}^{2}\zeta^{-\frac{1}{10}}.
\end{split}
\end{equation*}
According to  \eqref{wanfanhou} and Lemma \ref{kam} with
$\mathcal{P}^{-1} (e^{c_{0}-\delta_{1}c_{1}}+\delta_{1}^{2}P_{1})\mathcal{P} = \Delta  +\mathcal{P}^{-1}\delta_{1}^{2}P_{1}(\theta)\mathcal{P}$, we have
\begin{equation}
\begin{split}
&\ \ \ \ |\rho(\alpha,e^{c_{0}-\delta_{1}c_{1}}+\delta_{1}^{2}P_{1})-\sqrt{\det(c_{0}-\delta_{1}c_{1})}| \\
&= |\rho(\alpha,\Delta  +\mathcal{P}^{-1}\delta_{1}^{2}P_{1}(\theta)\mathcal{P} )-\rho(\alpha, \Delta)|\\
&\leq \delta_{1}^{2}\|\mathcal{P}\|^{2}\|P_{1}\|_{\hat{k}}\\
&\leq \zeta^{\frac{9}{5}}\times  4\|B\|_{k_{0}}^{2}\zeta^{-\frac{1}{10}} \times ( 53D_{\tau}^{2}\gamma^{-6}\|B\|_{k_{0}}^{4} +\zeta^{-\frac{9}{10}}\zeta^{2}\|B\|_{k_{0}}^{2} )\\
&\leq 240 D_{\tau}^{2}\gamma^{-6}\|B\|_{k_{0}}^{6}\zeta^{\frac{17}{10}}.
\end{split} \label{c2}
\end{equation}
By  (\ref{tiaojian}), we have
\begin{equation*}
240D_{\tau}^{2}\gamma^{-6}\|B\|_{k_{0}}^{6}\zeta^{\frac{7}{10}} \leq 1
\end{equation*}
then combine with (\ref{c1}) and  (\ref{c2}), we have
\begin{equation*}
\begin{split}
\rho(\alpha, e^{c_{0}-\delta_{1}c_{1}}+\delta_{1}^{2}P_{1})&\geq |\rho(\alpha,\Delta)|- |\rho(\alpha,\Delta  +\mathcal{P}^{-1}\delta_{1}^{2}P_{1}(\theta)\mathcal{P} )-\rho(\alpha, \Delta)|\\
& \geq \sqrt{3}\zeta -240 D_{\tau}^{2}\gamma^{-6}\|B\|_{k_{0}}^{6}\zeta^{\frac{17}{10}}\\
&\geq \sqrt{3}\zeta-\zeta >0.
\end{split}
\end{equation*}
recall that $\zeta \leq n^{-(k-62\tau)}_{J}$, we have
\begin{equation*}
|I_{J}(V)|\leq \delta_{1}=\zeta^{\frac{9}{10}} \leq n^{-\frac{9k}{10}+56\tau}_{J}, \ \ \ \forall  n_J\in\mathcal{K}.
\end{equation*}
This brings the proof of the upper bound estimate to its conclusion.

As we move forward, let's consider about the gap's lower bound estimate.
Let $\delta_2:=\zeta^{\frac{11}{10}}$. We will  demonstrate that $|I_{J}(V)| \geq \delta_2$.
Note that
$$\delta_2^2\left|[B^2_{11}][B^2_{12}]-[B_{11}B_{12}]^2\right|\leq  2 \zeta^{\frac{11}{5}}\|B\|_{k_0}^{4},$$
and by Lemma \ref{poly}, one has
$\delta_2[B^2_{11}]\zeta\geq \frac{1}{4} \zeta^{\frac{21}{10}} \|B\|_{k_0}^{-2}$.

As a result of \eqref{202307141} and \eqref{202307142}, we obtain
$\|B\|_{k_0}^{6}\zeta^{\frac{1}{10}}\leq \frac{1}{40}$. This implies
\begin{equation}\label{d2}
d(\delta_2)=-\delta_2[B^2_{11}]\zeta+\delta_2^2\left([B^2_{11}][B^2_{12}]-[B_{11}B_{12}]^2\right)<-\frac{1}{5}\zeta^{\frac{21}{10}} \|B\|_{k_0}^{-2},
\end{equation}
and hence
\begin{equation}\label{upper_determinant}
\sqrt{-d(\delta_2)}> \frac{1}{\sqrt{5}}\zeta^{\frac{21}{20}} \|B\|_{k_0}^{-1}.
\end{equation}

In view of Proposition 18 of \cite{Puig06}, there exists $\mathcal{P} \in \mathrm{SL}(2, \mathbb{R})$, with $\|\mathcal{P}\| \leq 2\left(\frac{\left|c_0-\delta_2 c_1\right|}{\sqrt{-d\left(\delta_2\right)}}\right)^{\frac{1}{2}}$ such that
$$
\mathcal{P}^{-1} e^{c_0-\delta_2 c_1} \mathcal{P}=\left(\begin{array}{cc}
e^{\sqrt{-d\left(\delta_2\right)}} & 0 \\
0 & e^{-\sqrt{-d\left(\delta_2\right)}}
\end{array}\right) .
$$
Since $\|B\|_{k_0}^{6}\zeta^{\frac{1}{10}}\leq \frac{1}{40}$, we have
$$
\left|c_0-\delta_2 c_1\right| \leq \zeta+\zeta^{\frac{11}{10}}(1+\zeta)\|B\|_{k_0}^{2} \leq 2 \zeta.
$$
Together with  $(\ref{upper_determinant})$, we have
$$\frac{|c_0-\delta_2 c_1|}{\sqrt{-d(\delta_2)}}\leq \frac{\sqrt{5}\cdot 2  \zeta}{\zeta^{\frac{21}{20}} \|B\|_{k_0}^{-1}}= 2\sqrt{5} \|B\|_{k_0}  \zeta^{-\frac{1}{20}}.$$
By Lemma \ref{kam}, we have
$$
\begin{aligned}
\mathcal{P}^{-1} \delta_2^2 P_1(\theta) \mathcal{P} & \leq 8 \sqrt{5}\|B\|_{k_0} \zeta^{-\frac{1}{20}} \zeta^{\frac{11}{5}}\left(53 D_\tau^2 \gamma^{-6}\|B(\theta)\|_{k_0}^4+\zeta^{-\frac{11}{10}} \zeta^2\|B\|_{k_0}^2\right) \\
& \leq 480 \sqrt{5} D_\tau^2 \gamma^{-6}\|B\|_{k_0}^5 \zeta^{\frac{43}{20}}
\leq -d(\delta_2),
\end{aligned}
$$
where the last inequality we use \eqref{tiaojian} and \eqref{d2}.
Consequently,  the cocycle $(\alpha, e^{c_0-\delta_2 c_1}+\delta_2^2 P_1)$ is uniformly hyperbolic, and $E_{n_J}^{+}-\delta_2 \not \in \Sigma_{V,\alpha}$, which means that $|I_{J}(V)| \geq \delta_2=\zeta^{\frac{11}{10}},$ which, together with \eqref{202307142}, yields
$$|I_{J}(V)| \geq n^{-\frac{11}{10}k-6\tau}_{J}.$$

\end{proof}
\section*{Acknowledgement}
J. He and H. Cheng were supported by NSFC grant (12001294). They would like to
give their thanks to Q. Zhou for useful discussions.

\bibliographystyle{plain}
\bibliography{202307cantorspectum}

\begin{thebibliography}{10}

\bibitem{Avilaa09}
A.~Avila.
\newblock On the spectrum and {L}yapunov exponent of limit periodic
  {S}chr\"{o}dinger operators.
\newblock {\em Comm. Math. Phys}, 288(3):907--918, 2009.

\bibitem{Artura10}
A.~Avila.
\newblock Almost reducibility and absolute continuity {I}.
\newblock {\em arXiv:1006.0704}, 2010.

\bibitem{AvilaA15}
A.~Avila.
\newblock Global theory of one-frequency {S}chr\"{o}dinger operators.
\newblock {\em Acta Math.}, 215(1):1--54, 2015.

\bibitem{Artura23}
A.~Avila.
\newblock K{AM}, {L}yapunov exponents and the spectral dichotomy for typical
  one-frequency {S}chr\"{o}dinger operators.
\newblock {\em arXiv:2307.11071}, 2020.

\bibitem{AvilaD08}
A.~Avila and D.~Damanik.
\newblock Absolute continuity of the integrated density of states for the
  almost {M}athieu operator with non-critical coupling.
\newblock {\em Invent. Math.}, 172(2):439--453, 2008.

\bibitem{Avila2009TheTM}
A.~Avila and S.~Jitomirskaya.
\newblock The ten martini problem.
\newblock {\em Ann. Math.}, 170:303--342, 2009.

\bibitem{Avila2003ReducibilityON}
A.~Avila and R.~Krikorian.
\newblock Reducibility or nonuniform hyperbolicity for quasiperiodic
  {S}chr\"odinger cocycles.
\newblock {\em Ann. Math.}, 164:911--940, 2003.

\bibitem{Avron81}
J.-E. Avron and B.~Simon.
\newblock Transient and recurrent spectrum.
\newblock {\em J. Funct. Anal.}, 43(1):1--31, 1981.

\bibitem{Cai2017ReducibilityOF}
A.~Cai and L.~Ge.
\newblock Reducibility of finitely differentiable quasi-periodic cocycles and
  its spectral applications.
\newblock {\em J. Dyn. Differ. Equ.}, 34:2079 -- 2104, 2017.

\bibitem{CAI2021109035}
A.~Cai and X.~Wang.
\newblock Polynomial decay of the gap length for {$C^k$} quasi-periodic
  {S}chr\"{o}dinger operators and spectral application.
\newblock {\em J. Funct. Anal.}, 281(3):Paper No. 109035, 30, 2021.

\bibitem{Chulaevsky1989AndersonLF}
V.~Chulaevsky and Y.~Sinai.
\newblock Anderson localization for the 1-d discrete {S}chr{\"o}dinger operator
  with two-frequency potential.
\newblock {\em Comm. Math. Phys}, 125:91--112, 1989.

\bibitem{Damanik2015}
J.~Fillman D.~Damanik and M.~Lukic.
\newblock Limit-periodic continuum {S}chr\"odinger operators with zero measure
  cantor spectrum.
\newblock {\em arXiv: Spectral Theory}, 2015.

\bibitem{Thouless8}
M.~Kohmoto D.~J.~Thoulessand, M.~P. Nightingale, and M.~den Nijs.
\newblock Quantized hall conductance in a two-dimensional periodic potential.
\newblock {\em Phys. Rev. Lett.}, 49:405--408, Aug 1982.

\bibitem{Dinaburg75}
E.~I. Dinaburg and Ja.~G. Sinai.
\newblock The one-dimensional {S}chr\"{o}dinger equation with quasiperiodic
  potential.
\newblock {\em Funkcional. Anal. i Prilo\v{z}en.}, 9(4):8--21, 1975.

\bibitem{Eliasson92}
L.~H. Eliasson.
\newblock Floquet solutions for the {$1$}-dimensional quasi-periodic
  {S}chr\"{o}dinger equation.
\newblock {\em Comm. Math. Phys}, 146(3):447--482, 1992.

\bibitem{Krikorian09}
B.~Fayad and R.~Krikorian.
\newblock Rigidity results for quasiperiodic {${\rm SL}(2,\Bbb R)$}-cocycles.
\newblock {\em J. Mod. Dyn.}, 3(4):497--510, 2009.

\bibitem{Goldstein2005OnRA}
M.~Goldstein and W.~Schlag.
\newblock On resonances and the formation of gaps in the spectrum of
  quasi-periodic {S}chr{\"o}dinger equations.
\newblock {\em Ann. Math.}, 173:337--475, 2011.

\bibitem{H83}
M.-R. Herman.
\newblock Une m\'ethode pour minorer les exposants de {L}yapounov et quelques
  exemples montrant le caract\`ere local d'un th\'eor\`eme d'{A}rnol'd et de
  {M}oser sur le tore de dimension {$2$}.
\newblock {\em Comment. Math. Helv.}, 58(3):453--502, 1983.

\bibitem{HouS19}
X.~Hou, Y.~Shan, and J.~You.
\newblock Construction of quasiperiodic {S}chr{\"o}dinger operators with cantor
  spectrum.
\newblock {\em Annales Henri Poincar{\'e}}, 20:3563 -- 3601, 2019.

\bibitem{HouY12}
X.~Hou and J.~You.
\newblock Almost reducibility and non-perturbative reducibility of
  quasi-periodic linear systems.
\newblock {\em Invent. Math.}, 190(1):209--260, 2012.

\bibitem{Bellissard1982}
D.~Bessis J.~Bellissard and P.~Moussa.
\newblock Chaotic states of almost periodic {S}chr{\"o}dinger operators.
\newblock {\em Phys. Rev. Lett.}, 49:701--704, 1982.

\bibitem{Johnson82}
R.~Johnson and J.~Moser.
\newblock The rotation number for almost periodic potentials.
\newblock {\em Comm. Math. Phys}, 84(3):403--438, 1982.

\bibitem{Klitzing80}
K.~V. Klitzing, G.~Dorda, and M.~Pepper.
\newblock New method for high-accuracy determination of the fine structure
  constant based on quantized hall resistance.
\newblock {\em Phys. Rev. Lett.}, 45(6):050501, 7, 1980.

\bibitem{Krikorian04}
R.~Krikorian.
\newblock Reducibility, differentiable rigidity and {L}yapunov exponents for
  quasi-periodic cocycles on $\mathbb{T}\times {SL}(2,\mathbb{R})$.
\newblock {\em Arxiv}, 2004.

\bibitem{zhaozhiyan}
M.~Leguil, J.~You, Z.~Zhao, and Q.~Zhou.
\newblock Asymptotics of spectral gaps of quasi-periodic {S}chr\"odinger
  operators.
\newblock {\em Arxiv}, 2017.

\bibitem{Moser84}
J.~Moser and J.~P\"{o}schel.
\newblock An extension of a result by {D}inaburg and {S}ina\u{\i} on
  quasiperiodic potentials.
\newblock {\em Comment. Math. Helv.}, 59(1):39--85, 1984.

\bibitem{Moser1981}
J.~K. Moser.
\newblock An example of a {S}chr\"odinger equation with almost periodic
  potential and nowhere dense spectrum.
\newblock {\em Comment. Math. Helv.}, 56:198--224, 1981.

\bibitem{puig}
J.~Puig.
\newblock Cantor spectrum for the almost {M}athieu operator.
\newblock {\em Comm. Math. Phys}, 244(2):297--309, 2004.

\bibitem{Puig06}
J.~Puig.
\newblock A nonperturbative {E}liasson's reducibility theorem.
\newblock {\em Nonlinearity}, 19(2):355--376, 2006.

\bibitem{SIMON1982463}
B.~Simon.
\newblock Almost periodic {S}chr{\"o}dinger operators: A review.
\newblock {\em Adv. Appl. Math.}, 3(4):463--490, 1982.

\bibitem{Simon1983KotaniTF}
B.~Simon.
\newblock Kotani theory for one dimensional stochastic jacobi matrices.
\newblock {\em Comm. Math. Phys}, 89:227--234, 1983.

\bibitem{Sodin95}
M.~Sodin and P.~Yuditskii.
\newblock Almost periodic {S}turm-{L}iouville operators with {C}antor
  homogeneous spectrum.
\newblock {\em Comment. Math. Helv.}, 70(4):639--658, 1995.

\bibitem{Sodin97}
M.~Sodin and P.~Yuditskii.
\newblock Almost periodic {J}acobi matrices with homogeneous spectrum,
  infinite-dimensional {J}acobi inversion, and {H}ardy spaces of
  character-automorphic functions.
\newblock {\em J. Geom. Anal.}, 7(3):387--435, 1997.

\bibitem{Wang2014CantorSF}
Y.~Wang and Z.~Zhang.
\newblock Cantor spectrum for a class of ${C}^2$ quasiperiodic {S}chr\"odinger
  operators.
\newblock {\em arXiv: Dynamical Systems}, 2014.

\bibitem{You2019}
J.~You.
\newblock Quantitative almost reducibility and its applications.
\newblock In {\em Proceedings of the International Congress of Mathematicians
  (ICM 2018) (In 4 Volumes) Proceedings of the International Congress of
  Mathematicians 2018}, pages 2113--2135. World Scientific, 2018.

\end{thebibliography}

\end{document}